\documentclass[12pt]{article}
\usepackage[paper=letterpaper,margin=0.9in]{geometry}

\usepackage{tocloft}
\usepackage{amsmath,amssymb,mathrsfs,tensor,bbold,upgreek,cases}
\usepackage{amsthm,amsfonts,mathtools,bm,slashed}
\usepackage{listings}
\usepackage{float,graphicx,comment,subfig}
\usepackage[hidelinks]{hyperref}
\usepackage{cite}
\usepackage[font=small,labelfont=bf]{caption}

\numberwithin{equation}{section}

\usepackage{afterpage,placeins}
\usepackage[dvipsnames]{xcolor}

\usepackage{listings}
\usepackage[lined,ruled]{algorithm2e}

\usepackage{tikz-cd}
\tikzcdset{arrow style=tikz, diagrams={>=stealth}}
\usepackage{tikz}
\usepackage{pgfplots}
\usepgfplotslibrary{fillbetween}
\usetikzlibrary{arrows,shapes,positioning}
\usetikzlibrary{decorations.markings}
\usetikzlibrary{calc}
\definecolor{error1}{rgb}{1.0, 0.25, 0.25}
\definecolor{quad4}{RGB}{30,171,154}
\definecolor{quad1}{RGB}{222,0,109}
\definecolor{quad3}{RGB}{175,209,48}
\definecolor{quad2}{RGB}{235,67,30}

\def\be{\begin{equation}}
\def\ee{\end{equation}}
\def\bea{\begin{align}}
\def\eea{\end{align}}

\DeclareMathOperator*{\argmin}{arg\,min}

\DeclareMathOperator{\dom}{dom}

\DeclareMathOperator{\zer}{zer}

\DeclareMathOperator{\prox}{prox}

\newcommand{\Order}{\mathcal{O}}
\newcommand{\Mobs}{M_{\textnormal{obs}}}
\newcommand{\PP}{\mathcal{P}_{\Omega}}

\newcommand{\bigO}{\mathcal{O}}
\newcommand\smallO{
  \mathchoice
    {{\scriptstyle\mathcal{O}}}% \displaystyle
    {{\scriptstyle\mathcal{O}}}% \textstyle
    {{\scriptscriptstyle\mathcal{O}}}% \scriptstyle
    {\scalebox{.7}{$\scriptscriptstyle\mathcal{O}$}}%\scriptscriptstyle
  }
\newcommand{\F}{\varphi}

\newcommand{\E}{\mathbb{E}}

\newcommand{\B}{\mathcal{B}}
\theoremstyle{remark}

\newcommand\mailto[1]{\href{mailto:#1}{\tt #1}}

\makeatletter
\def\blfootnote{\gdef\@thefnmark{}\@footnotetext}
\makeatother

\begin{document}
\thispagestyle{empty}
\renewcommand{\thefootnote}{\fnsymbol{footnote}}

\vspace*{0.5cm}

\begin{center}

{\bf{\large Gradient flows and proximal splitting methods:
A unified view on  \\[.5em] accelerated and stochastic optimization}}

\vspace{1.0cm}

{\bf Guilherme Fran\c ca,}$^{\!1,2,}%
\footnote{\mailto{guifranca@gmail.com}}$%
{\bf ~Daniel P. Robinson,}$^{\!3}$%
{\bf ~and Ren\' e Vidal}$^{\,2}$

\vspace{0.5cm}

{\it
${}^{1}$University of California, Berkeley, California 94720, USA\\[.2em]
${}^{2}$Johns Hopkins University, Maryland 21218, USA\\[.2em]
${}^{3}$Lehigh University, Pennsylvania 18015, USA
}

\vspace{1.5cm}

{\bf Abstract}
\vspace{-1em}
\end{center}

Optimization is at the heart of machine learning, statistics  and many
applied scientific disciplines. It also has a long history in physics,
ranging from the minimal action principle to finding  ground states of
disordered systems such as spin glasses.
Proximal algorithms form a class of  methods that are broadly applicable
and  are particularly well-suited  to nonsmooth, constrained,
large-scale, and distributed optimization problems.
There are essentially five proximal algorithms currently known, each
proposed in seminal work:
Forward-backward splitting,
Tseng splitting,
Douglas-Rachford,
alternating direction method of multipliers,
and the more recent Davis-Yin.
These methods sit on a higher level of abstraction compared to
gradient-based ones, with
deep roots in nonlinear functional analysis.
In this paper we show that all of these methods
are actually different
discretizations of a single differential equation,
namely, the simple \emph{gradient flow} which dates back to Cauchy (1847).
An important aspect behind
many of the success stories in machine learning relies on ``accelerating''
the convergence of first-order methods. However, accelerated
methods are  notoriously difficult to analyze, counterintuitive,
and without
an underlying guiding principle.
We show that similar discretization schemes applied to
Newton's equation with an additional dissipative force,
which we refer to as
\emph{accelerated gradient flow}, allow us to obtain accelerated
variants of all these proximal algorithms---the majority of which are
new
although some recover known cases in the literature.
Furthermore, we extend these methods to stochastic settings,
allowing  us to make connections with Langevin and Fokker-Planck equations.
Similar ideas apply to
gradient descent, heavy ball, and Nesterov's method which are simpler.
Our results therefore provide a unified framework from which
several important optimization methods are nothing but \emph{simulations} of
classical dissipative systems.

\renewcommand*{\thefootnote}{\arabic{footnote}}
\setcounter{footnote}{0}
\setcounter{page}{0}
\newpage
\thispagestyle{empty}
\tableofcontents
\newpage

\section{Introduction}

The simplest algorithm to solve a \emph{smooth} optimization  problem
\begin{equation} \label{optimization}
\min_{x \in \mathbb{R}^n} \F(x)
\end{equation}
dates back to Cauchy \cite{Cauchy:1847}. It is the well-known
\emph{gradient descent} method, $x_{k+1} = x_k - h \nabla \F(x_k)$,
where $h>0$ is the step size and $k=0,1\dotsc$ is the iteration number.
Clearly, gradient descent corresponds to
an \emph{explicit} Euler discretization of the
\emph{gradient flow}:
\begin{equation} \label{gradflow}
\dot{x} = - \nabla \F(x),
\end{equation}
where $x \equiv x(t)$.
To minimize  \emph{nonsmooth} and \emph{composite} functions,
a series of milestone papers introduced algorithms based on
\emph{proximal operators}, which do not require explicit gradient computations.
For instance,
the \emph{Douglas-Rachford} algorithm \cite{Douglas:1956} was proposed
in the 50's to solve the heat equation but nowadays is a standard
optimization
method with important applications.
Closely related is the \emph{alternating direction method of multipliers}
(ADMM), introduced independently by
Glowinsky and Marrocco \cite{Glowinsky:1975}
and Gabay and Mercier \cite{Gabay:1976} in the 70's---ADMM
has been gaining significant interest in machine learning
during the last decade \cite{Boyd:2011}.
Another method---that plays an important role in signal processing---is
the \emph{forward-backward splitting} introduced by
Lions and Mercier \cite{Lions:1979} also in the 70's.
These were the only known proximal-based methods for almost 30 years,
until  Tseng \cite{Tseng:2000} proposed a slight
modification of the
latter known as \emph{forward-backward-forward splitting}.
Such methods are designed to minimize composite functions
$\F(x) = \F_1(x) + \F_2(x)$, where both $\F_1, \F_2$ are allowed
to be nonsmooth.
Finding an algorithm that minimizes
$\F(x) = \F_1(x) + \F_2(x) + \F_3(x)$, where
only $\F_3$ is assumed to be smooth, and which cannot
be reduced to any of the previously known methods, was a longstanding
problem that has
been recently solved by Davis and Yin \cite{Davis:2017}.
These five  algorithms compose the list of fundamental
proximal algorithms currently known---many other methods
are variations of these basic themes. Such proximal methods
can be derived from operator splitting techniques \cite{Ryu:2016,Parikh:2013},
which have origins in the works of Browder
\cite{Browder:1963,Browder:1963b,Browder:1963c} and
Minty \cite{Minty:1962}, although
nowadays form an entire field of research in convex analysis, optimization,
and nonlinear functional analysis \cite{BauschkeBook,Zeidler2}.

Perhaps surprisingly, we provide a simple yet unified
perspective on these distinct methods: All of them are
different discretizations of the simple gradient flow \eqref{gradflow}.
More precisely, they are  first order integrators that preserve
the critical points of this ODE.

``Acceleration strategies'' in  the context of optimization have proved to be
powerful and are behind many of the empirical success stories in  machine
learning, such as in training neural networks.
The basis of accelerated \emph{gradient-based} methods can be traced back to
Polyak \cite{Polyak:1964} and Nesterov~\cite{Nesterov:1983}---both can be
seen as accelerated versions of gradient descent. Although neither are
intuitive in their precise design, it has recently been shown that both
can be obtained as \emph{explicit} discretizations of a
second-order ODE
\cite{Candes:2016,Wibisono:2016,Franca:2019}.
This continuous-time perspective on optimization methods
is quite recent and has
helped demystify the ``magic" of acceleration techniques.
However, the construction and analysis of accelerated methods is still
obscure, without an underlying guiding principle---e.g., it is not clear how to
``accelerate'' a given algorithm. \emph{Accelerated proximal-based} methods
are even less known, and can play an important role since they may
enjoy improved stability and be applicable in more general situations.
Moreover, from a mathematical standpoint, such methods  have  interesting
connections with nonlinear functional analysis.  As we will show,
several accelerated variants of each of the above mentioned
proximal algorithms
can be obtained as different discretizations of
\begin{equation}
\label{secode}
\ddot{x} + \eta(t)\dot{x} = -\nabla \F(x) ,
\end{equation}
with a suitable choice of damping function $\eta(t) > 0$.
The resulting methods---most of which are new in the
literature---are  first-order integrators that preserve critical points
of  this ODE.
Therefore, the above \emph{classical dissipative system}
has a fundamental importance to
optimization.
Note that this is nothing but
Newton's equation with an additional dissipative force, $-\eta(t) \dot{x}$.
When $\eta(t) = \eta$ is constant and $\F = \omega^2 x^2 / 2$
it reduces to the Caldirola-Kanai  oscillator \cite{Caldirola:1941,Kanai:1948},
which is the classical limit of the seminal Caldeira-Legget
model \cite{Caldeira:1981}.

Our approach makes connections between optimization and
\emph{splitting methods} for ODEs \cite{McLachlan:2002}.
Interestingly, ADMM and its accelerated variants arise as a
\emph{rebalanced splitting}, which is a recent technique
designed to preserve critical points
\cite{Speth:2013}---the so-called dual variable of ADMM, originally introduced
as a Lagrange multiplier, is precisely the balance coefficient
of Ref. \cite{Speth:2013}. The other methods we consider also preserve
critical points,   but for different reasons, which
in turn sheds light on the connections between ODE
splitting  and operator  splitting  ideas from convex  analysis.

\emph{Stochastic optimization} is an important ingredient in the
machine learning
toolbox to reduce
the computational burden in training high-dimensional models
over large datasets.
By introducing stochastic gradients or
stochastic  proximal operators into these methods, instead
of systems \eqref{gradflow} and \eqref{secode}, their
continuum limit become an \emph{overdamped} or
\emph{underdamped} Langevin equation, respectively. The probability
distribution of such stochastic processes is
described by a Fokker-Planck equation.
Therefore, there is a close connection between deterministic optimization
and dissipative classical mechanics, as well as
stochastic optimization and nonequilibrium statistical mechanics.

This paper is organized as follows.
In Sec.~\ref{sec:resolvent}, we introduce basic concepts related
to proximal operators---or monotone operators and
their regularizations more generally---and illustrate how they naturally arise
from \emph{implicit} discretizations of ODEs.
In Sec.~\ref{sec:rates}, we show relevant details about the
dynamics of the gradient flows \eqref{gradflow} and \eqref{secode}.
In Sec.~\ref{sec:admm}, we introduce a slight variation of
the balanced and rebalanced splitting schemes \cite{Speth:2013}
to then show how---accelerated---ADMM arises from this approach.
In Sec.~\ref{sec:davis_yin}, we derive extensions of
Davis-Yin, which is known to generalize both
forward-backward and Douglas-Rachford.
%hence our results immediately imply connections
%for these methods as well.
%
In Sec.~\ref{sec:tseng},
we introduce accelerated extensions of Tseng's splitting to complete the list.
The focus of this paper is on discretizations of the second-order
gradient  flow  \eqref{secode} since this allows us to construct entire
new families of accelerated methods that
generalize the existing  ones.
Moreover, the known methods that are related to
\eqref{gradflow} follow as particular cases---more precisely, through
a high friction limit.
At this stage, we  briefly
summarize and interpret
these results from a physics perspective in Sec.~\ref{sec:classview}.
Then, in Sec.~\ref{sec:stochastic},
we shift gears and
show how one can extend---quite easily---these proximal-based methods  to
stochastic optimization  settings.
As a consequence, the connections with the continuous-time  formalism
are  promoted to SDEs of the Langevin type,
whose probability distribution are governed  by Fokker-Planck equations.
For completeness, in the Appendix
we show that stochastic gradient descent,
heavy ball, and Nesterov fit our general framework---such gradient-based methods
find widespread applications in machine learning but
they are actually simpler.
We provide numerical experiments in Sec.~\ref{sec:numerical}  that
support our  theoretical  conclusions,  and also illustrate  the  faster
convergence attained  by  the accelerated methods.
%in
%order to illustrate the speedup achieved by our proposed accelerated methods.
%
Our final remarks and a potential implications of
our analysis are presented in Sec.~\ref{sec:conclusion}.

\section[Resolvent, Yosida regularization, and proximal operator]{Resolvent,
Yosida regularization,  \\ and proximal operator}
\label{sec:resolvent}

We start by introducing fundamental concepts from nonlinear
functional analysis \cite{BauschkeBook,Zeidler2} since this is the formalism
in which proximal algorithms are generally discussed.
We avoid excessive formalism in the paper, but here we
give a roadmap to further abstract our  analysis.

The \emph{resolvent} of an operator $A$ can be defined as
\begin{equation}
\label{resolvent}
J_{\lambda A} \equiv \left( I + \lambda A \right)^{-1} ,
\end{equation}
where $\lambda$ is the so-called spectral parameter.
Even though $\lambda$ can be complex,
we will only need $\lambda \in \mathbb{R}$.
Another useful concept is the \emph{Yosida regularization} of
$A$:
\begin{equation}
  \label{yosida}
  A_\lambda \equiv \lambda^{-1} (I - J_{\lambda A}).
\end{equation}
Let $H$ be a Hilbert space with inner product $\langle \cdot | \cdot \rangle:
H \times H \to \mathbb{C}$. A \emph{multivalued}
map $A: H \rightrightarrows H$, with
$\dom A \equiv \{ x \in H \, \vert \, Ax \ne \emptyset \}$, is said to
be \emph{monotone} if and only if
\begin{equation} \label{mon}
\langle A y - Ax | y - x\rangle \ge 0 \qquad \mbox{ for all
$x,y\in\dom A$}.
\end{equation}
A monotone operator is said to be \emph{maximal} if
no enlargement of its graph is possible. It turns out that every monotone
operator admits a maximal extension, thus we henceforth assume
that all operators
are maximal monotone.
What matters for us is that in this case
the resolvent \eqref{resolvent} is
single-valued, i.e., $J_{\lambda A}: H \to H$ is a function.
Moreover, $x^\star$ is a zero of $A$,
namely $x^\star \in \zer(A) \equiv \{x\in H \, \vert \, 0\in Ax \}$,
if and only if  $J_{\lambda A}(x^\star) = x^\star$.
Thus, zeros of $A$ are
fixed points of the resolvent, $J_{\lambda A}$, and vice-versa.
Consequently, the Yosida regularization \eqref{yosida} is
also single-valued, and $x^\star$ is a zero of $A$ if and only
if $A_\lambda x^\star = 0$, so that
zeros of $A$ are also zeros of $A_\lambda$ and vice-versa.
The advantage of working with the Yosida regularization is that it allows us
to deal with multivalued operators by considering single-valued operators.
Indeed, it can be shown
that $A_{\lambda} x \to A_0 x$
as $\lambda \downarrow 0$,
where $A_0 x$ is the element
of minimal norm in the set $A x$.

These ideas can be made more intuitive by considering a function
$\F : \mathbb{R}^n \to \mathbb{R}$, which for the moment we  assume to be
differentiable.
The function $\F$ is convex if and only if its gradient,
$A = \nabla \F$, is (maximal)  monotone.
In this case, the resolvent \eqref{resolvent}
becomes the so-called \emph{proximal operator}:%
\footnote{This can be easily seen as follows. The solution  $y$ of
\eqref{prox} obeys $\nabla \F(y) + (1/\lambda)(y - x) = 0$, i.e.,
$(I + \lambda  \nabla \F)y = x$, which from Eq. \eqref{resolvent} gives
$y = J_{\lambda \nabla \F} x$. We replaced $\nabla \F \mapsto \partial \F$ in
\eqref{prox} anticipating generalization to nonsmooth cases.}
\begin{equation}
\label{prox}
\begin{split}
J_{\lambda \partial \F}(x) &\equiv
\prox_{\lambda \F}(x)   \\ &\equiv \argmin_y \left(
\F(y) + \dfrac{1}{2\lambda}\| y - x\|^2\right).
\end{split}
\end{equation}
Then Eq. \eqref{yosida} becomes
$(\nabla \F)_\lambda(x) = \lambda^{-1}\big[ x-\prox_{\lambda \F}(x) \big]$, which
is the gradient of the \emph{Moreau envelope} $\F_\lambda$,
i.e., $(\nabla \F)_\lambda(x) = \nabla \F_\lambda(x)$
where
\begin{equation}\label{Moreau}
\F_{\lambda}(x) \equiv \min_{y}\left(\F(y) +
\dfrac{1}{2\lambda}\| y - x\|^2 \right).
\end{equation}
When $\F$ is \emph{nonsmooth} its gradient $\nabla \F$ is ill-defined.
However, there exists a generalization which is the notion of
subdifferential set;
it is defined as
$\partial \F(x) \equiv \{ g \in \mathbb{R}^n \, | \, \F(y) - \F(x)
\ge \langle g | y - x \rangle \, \mbox{for all $y \in \mathbb{R}^n$}\}$.
If $\F$ is differentiable, then $\partial \F(x) = \{ \nabla \F(x) \}$.
We thus see that even though $\F$ may not be differentiable,
its Moreau envelope  always is, and we can thus treat the problem with
standard calculus on $\F_\lambda$.  We have
$\lim_{\lambda \downarrow 0}\nabla \F_\lambda(x) \in \partial \phi(x)$,
and this limit is the vector of minimal norm in the subdifferential
set $\partial \F(x)$.

More generally, our results in this paper show that all the
previously mentioned algorithms correspond to discretizations
of the \emph{differential inclusion} \cite{Zeidler2}
\begin{equation} \label{inc1}
\dot{x} \in - A x ,
\end{equation}
for a monotone operator $A : H \rightrightarrows H$ that is
composite, $A = A_1 + A_2 + A_3$.
Similarly, the accelerated variants of these algorithms
are related to the second-order differential inclusion
\begin{equation} \label{inc2}
\ddot{x} + \eta(t)\dot{x} \in - A x.
\end{equation}
However, dealing with differential inclusions, i.e., nonsmooth dynamical
systems, involve several technicalities.
A simple way to avoid the  issue is to focus on their Yosida regularizations,
namely
\begin{equation}
\dot{x} = - A_\mu x \label{yos1}
\end{equation}
and
\begin{equation}
\ddot{x} + \eta(t)\dot{x} = - A_\mu x, \label{yos2}
\end{equation}
respectively,  which are well-posed ODEs.
(Note that $x(t)$ depends on  $\mu > 0$, which we omit
in the notation.)
At the end of the day one can  take the limit $\mu \downarrow 0$
to recover results for \eqref{inc1} and \eqref{inc2} \cite{Zeidler2}.
In the context of nonsmooth optimization this means considering
the gradient of the Moreau envelope, $\nabla \F_\mu$, instead of
the subdifferential, $\partial \F$---this point will be further
clarified below.

As a warmup, and also to introduce the basic building blocks of our approach,
let us show a simple example on how to derive a proximal algorithm
from \eqref{inc1}, or equivalently \eqref{yos1}.
Consider an \emph{implicit} discretization of \eqref{inc1}:
\begin{equation}
\dfrac{x_{k+1} - x_k}{h} + \bigO(h) \in - A x_{k+1}.
\end{equation}
Using the resolvent \eqref{resolvent}, and neglecting $\bigO(h^2)$ terms,
we can solve this recurrence as
\begin{equation} \label{proxpointunreg}
x_{k+1} = J_{h A} x_k.
\end{equation}
This algorithm finds zeros of the monotone operator $A$.
For a nonsmooth function $\F$,  we set
$A  = \partial \F$ to obtain
\begin{equation}  \label{proximal_point}
x_{k+1} = \prox_{h \F}(x_k).
\end{equation}
This is the well-known  \emph{proximal gradient} method,
extensively studied in convex analysis and
optimization literatures.
Now, consider instead an analogous discretization of
the regularized ODE \eqref{yos1}:
\begin{equation}
\dfrac{x_{k+1} - x_k}{h} + \bigO(h) = - A_\mu x_{k+1}.
\end{equation}
From the resolvent \eqref{resolvent} we get
\begin{equation} \label{proxpointreg}
x_{k+1} = J_{h A_\mu}x_k .
\end{equation}
Employing the useful formula \cite{BauschkeBook}
\begin{equation}
\label{resolvent3}
J_{\lambda  A_\mu} =
(\lambda + \mu)^{-1}\left( \mu I +  \lambda J_{(\lambda + \mu)A} \right),
\end{equation}
we conclude that $J_{\lambda A_\mu} \to J_{\lambda A}$ as $\mu \downarrow 0$.
Thus \eqref{proxpointreg} becomes
precisely \eqref{proxpointunreg} in the limit $\mu \downarrow 0$.
In the case of a nonsmooth function $\F$, the RHS of \eqref{yos1} is
simply $-\nabla \F_\mu$, i.e., the gradient of Moreau envelope, whereas the
RHS of \eqref{inc1} is $-\partial \F$, i.e., the subdifferential set.
Hence, we could have considered a discretization of the gradient flow
\eqref{gradflow} with $\F \mapsto \F_\mu$ and then take the limit
$\mu \downarrow 0$. Even simpler, we could have actually ignored
nonsmoothness issues altogether and simply discretized
\eqref{gradflow}, replacing $\nabla \F \mapsto \partial \F$
where appropriate---this results in
updates in terms of the proximal operator \eqref{prox} and the procedure
can be formally justified by the above steps.

Next, let us consider a similar discretization approach
but for the second-order ODE
\eqref{yos2}. The differential operator on the LHS can be discretized as
\begin{equation}
\dfrac{x_{k+1} -2 x_k + x_{k-1}}{h^2} +
\eta_k \dfrac{x_k - x_{k-1}}{h} + \bigO(h).
\end{equation}
Defining
\begin{equation}\label{xhat}
\hat{x}_k \equiv x_k + \gamma_k(x_k - x_{k-1}), \qquad
\gamma_k \equiv (1 - h \eta_k),
\end{equation}
where $\eta_k \equiv \eta(t_k)$ is the discretized damping coefficient,
assumed to be only a function of time,
we obtain
\begin{equation} \label{eulersec}
\ddot{x}(t_k) + \eta(t_k)\dot{x}(t_k) =
\dfrac{x_{k+1} - \hat{x}_k}{h^2} + \bigO(h).
\end{equation}
This relation will prove extremely convenient
in pretty much all discretizations considered in this paper.
Note that it allows us to discretize the second-order system
\eqref{yos2} in very similar way as
the first-order system \eqref{yos1}---essentially, it suffices to replace
$x_k \mapsto \hat{x}_k$ and $h \mapsto h^2$.
Thus,  an implicit discretization
of \eqref{yos2}
yields $(x_{k+1} - \hat{x}_k)/h^2 = -A_\mu x_{k+1}$, which can be readily
solved with
the resolvent \eqref{resolvent} to obtain
\begin{equation} \label{aproxpointreg}
x_{k+1} = J_{h^2 A_\mu} \hat{x}_k.
\end{equation}
Taking $\mu \downarrow 0$ yields $x_{k+1} = J_{h^2 A}\hat{x}_k$.
For a nonsmooth function $\F$, we set
$A = \partial \F$ and replace
the resolvent by the proximal operator, hence obtaining
\begin{subequations} \label{a_proximal_point}
\begin{align}
x_{k+1} &= \prox_{h \F}(\hat{x}_k) \\
\hat{x}_{k+1} &= x_{k+1} + \gamma_{k+1}(x_{k+1}-x_{k}).
\label{a_proximal_point.u2}
\end{align}
\end{subequations}
Note that in  this last passage we redefined the step size,
$h^2 \mapsto h$,
and this should be also reflected in \eqref{xhat}.
The above algorithm corresponds to an ``accelerated version'' of the
proximal gradient
method \eqref{proximal_point}---in the next section we will show that
this second-order system
may indeed have faster convergence compared to
the first-order gradient flow.
Note also that the nonaccelerated method
\eqref{proximal_point} can
be  recovered from \eqref{a_proximal_point} by setting $\gamma_k = 0$.
Physically, this corresponds to a ``large friction limit'' as we will explain
in more detail shortly.
%%
%%%%%%%%%%%%%%
%\footnote{In this case $\hat{x}_k = x_k$; see Eq. \eqref{xhat}.
%Dropping $k$ for simplicity, note that
%$\eta = (1-\gamma)/h \to 1/h$ as $\gamma \to 0$. Thus
%$\eta \to \infty$ for sufficiently small step sizes $h \to 0$, which
%is a ``large friction limit.''
%Indeed, replacing $\eta \to 1/h$ into
%\eqref{eulersec} yields
%$\ddot{x} + \dot{x}/h = (x_{k+1} - x_k)/h^2 + \bigO(h)$ so that
%$\ddot{x}$ becomes negligible, i.e. this relation becomes
%$\dot{x} = (x_{k+1} - x_k)/h + \bigO(h)$ which is precisely
%a discretization of the first order system \eqref{yos1}.
%In addition, notice that \eqref{yos2} is Newton's equation in natural units
%where the mass $m=1$. Restoring the mass we have
%$m \ddot{x} + m \eta \dot{x} = - A_\mu x$. The large $\eta$ limit, often
%referred to as ``overdamped'' limit, can
%be obtained with $\eta \to 1 / m$ and taking $m\to 0$,
%thus yielding $\dot{x} = - A_\mu x$. In this sense
%\eqref{yos1} can be seen as a large friction
%limit---or zero-mass limit---of \eqref{yos2}; in this regime a very light
%particle has negligible acceleration due a large frictional force.
%\label{foot_high_friction}}
%%%%%%%%%%%%%%%%%

Already at this stage we have several new methods
encoded in \eqref{a_proximal_point}
due to the possibility of choosing different damping functions $\eta(t)$.
Reasonable choices (that will be justified in the next section) are a
\emph{constant damping},
\begin{equation} \label{hbdamp}
\eta(t) = \eta  \quad \implies \quad \gamma_k  = 1 - \sqrt{h} \, \eta,
\end{equation}
which is originally related to Polyak's heavy ball method \cite{Polyak:1964},
and a \emph{decaying damping},
\begin{equation} \label{nagdamp}
\eta(t) = r/t \quad (r\ge3) \quad \implies \quad \gamma_k = k/(k+r) ,
\end{equation}
where usually $r=3$ and this is originally
related to Nesterov's method \cite{Nesterov:1983,Candes:2016}.%
%%%%%%%%%%%%%%%%%%%%%%%
\footnote{The particular choice of $\gamma_k$ in \eqref{nagdamp} is to
maintain  consistency
with the optimization literature but $\gamma_k = k/(k+r) \approx 1 - r/k
= 1 - h \eta_k$ for
large $k$ which can equally be used.}
%%%%%%%%%%%%%%%%%%%%%%
However, other choices are possible and we allow
an arbitrary $\eta(t)$ in our discretizations.

\subsection{Large friction limit} \label{sec:high_friction}

We noted that setting $\gamma_k=0$ into \eqref{a_proximal_point}
yields the proximal gradient method \eqref{proximal_point},
which is a discretization of the first-order gradient
flow \eqref{gradflow}. We also said that this corresponds to
a ``large friction limit'' of the second-order gradient flow \eqref{secode}.
Since this idea will reappear later on, we provide more
details.

First, from a discrete-time viewpoint,
with $\gamma_k = 0$ we have $\hat{x}_k = x_k$.
Dropping  $k$ for simplicity, note that
the damping becomes
$\eta = (1-\gamma)/h \to 1/h$ as $\gamma \to 0$,  i.e.,
$\eta \to \infty$ as $h \to 0$. Therefore,
this corresponds indeed  to an ``overdamped limit.''
Moreover, replacing $\eta \to 1/h$ into
\eqref{eulersec} yields
\begin{equation}
h \ddot{x} + \dot{x} = \dfrac{x_{k+1} - x_k}{h} + \bigO(h^2).
\end{equation}
Note that $\ddot{x}$ becomes negligible for sufficiently small step sizes,
i.e.,
the above reduces to $\dot{x} = (x_{k+1} - x_k)/h + \bigO(h)$, which is precisely
a discretization of the LHS of Eq. \eqref{gradflow}---or
of Eq. \eqref{yos1} more generally.

Second, from a continuous-time viewpoint, note that Eq. \eqref{yos2} is
in natural units where the mass $m=1$. Restoring the mass
we have
\begin{equation} \label{foot_high_friction}
m \ddot{x} + m \eta \dot{x} = - A_\mu x.
\end{equation}
The overdamped limit can
be obtained with $\eta \to 1 / m$ and $m\to 0$, in which case
we recover Eq. \eqref{yos1}. It is in this
sense that the gradient flow is
a large friction limit---or a zero-mass limit---of the accelerated
gradient flow.
Intuitively, a very light particle has negligible acceleration due a
large frictional force, and the first-order system \eqref{yos1}
approximates the dynamics of the second-order system \eqref{yos2}
when $\eta$ is sufficiently large.
This same idea often appears in the theory  of
stochastic processes where a second-order
Langevin equation is well-approximated by a first-order Langevin equation
when the ``fluid'' representing  the heat bath is highly viscous.

\subsection{A note on nonsmoothness}
Above, we discretized the regularized ODEs \eqref{yos1} and
\eqref{yos2} and then  took the limit $\mu \downarrow 0$
to reduce the fixed point iterations to the case of
monotone operators. By choosing $A = \partial \F$ these
algorithms are appropriate for minimizing a nonsmooth function $\F$, through
the proximal operator \eqref{prox}.  Although we were careful in
taking ``nonsmoothness'' into account,
apart from
this $\mu$ limit, the discretization procedure is
exactly the same as if we had considered \eqref{gradflow} and \eqref{secode}
for which $\F$ is assumed to be differentiable. In other words,
everything works fine if we replace
$\nabla \F \mapsto \partial \F$ where appropriate---this point
was also noted in the discussion after Eq. \eqref{resolvent3}.
Moreover, even when
we split the operators, as we will do in the following, it is still possible
to introduce some parameter $\mu$ that justifies the procedure.
Therefore, to avoid unnecessary formalism,
hereafter we assume that $\F$  is differentiable for all practical purposes---%
one should keep in mind that this assumption can be removed by introducing
Yosida regularizations or Moreau envelopes and taking
the $\mu  \downarrow 0$ limit.

\section{Continuum dynamics}
\label{sec:rates}

Here we provide some details about the dynamics
of the gradient flows \eqref{gradflow} and \eqref{secode}.
Note that the first-order system \eqref{gradflow} yields the simplest
dynamics  that follows a descent direction on $\F$,
thus it is naturally suited for optimization purposes.
The second-order system \eqref{secode}
corresponds to its accelerated version in a classical
mechanical sense, and also follows a descent direction but can oscillate.
This is an actual dissipative system with Lagrangian
\begin{equation}
\mathcal{L} = \dfrac{1}{2} e^{{\theta}(t)} \| \dot{x}\|^2  -
e^{\theta(t)}\F(x),
\end{equation}
where $\dot{\theta}(t) \equiv \eta(t)$,
or equivalently with the explicit time-dependent Hamiltonian
\begin{equation}
\mathcal{H} = \dfrac{1}{2} e^{-{\theta}(t)} \| p\|^2  + e^{\theta(t)}\F(x) .
\end{equation}
The physical energy is given by $\mathcal{E} =
\tfrac{1}{2} \| \dot{x} \|^2 + \F(x)$
and dissipates at  a rate
$\dot{\mathcal{E}} = -\eta(t) \| \dot{x} \|^2 \le 0$, i.e.,
it decreases  monotonically with time, thus forcing trajectories to approach
the ground state $\F^\star \equiv \min_x \F(x)$.
In fact, $\mathcal{E}$ is a \emph{Lyapunov function}, enabling us to conclude
that the system is \emph{stable}\footnote{Stability means
that the trajectories stay nearby $x^\star$ for all times.}
on a minimizer  $x^\star \equiv \argmin_x \F(x)$.
This holds for any bounded damping $\eta(t) > 0$. In addition,
if $\eta(t) = \eta$ is constant then the system is
\emph{asymptotically stable}\footnote{This is stronger,
i.e., $x(t) \to x^\star$ as $t \to \infty$.
This result can be derived from LaSalle's invariance principle.}
around $x^\star$.

%We are working on natural units where
%$m=1$. Restoring the mass:
%$m \ddot{x} + m \eta \dot{x} = -\nabla \F(x)$.
%Its high friction limit can be obtained by
%redefining $\eta \to \eta / m$  and taking  $m \to 0$, which leads
%to $ \eta \dot{x} = - \nabla \F(x)$. In this sense, the gradient flow
%\eqref{gradflow} corresponds to a high friction limit of the
%accelerated gradient flow \eqref{secode}.
As explained in Sec.~\ref{sec:high_friction},
the gradient flow \eqref{gradflow}
corresponds to a large friction limit of the accelerated gradient
flow \eqref{secode}.
It is straightforward to show that the gradient flow is
\emph{asymptotically stable}
on $x^\star$.%
%%%%%%%%%%%%%%%%%%%%%%%%%%
\footnote{Indeed, consider the Lyapunov function
$\mathcal{E} \equiv \F(x) - \F^\star$.
We have $\mathcal{E} \ge 0$ and
$\dot{\mathcal{E}} = - \| \nabla \F\|^2 \le 0$, and note that such
inequalities are strict off a critical point.}
%%%%%%%%%%%%%%%%%%%%%%%%%
Thus, its trajectories
converge to $x^\star$, and so do trajectories of
\eqref{secode} with a constant damping.
However, when $\eta(t)$
is a decreasing function of time, such as in \eqref{nagdamp}, we can
only conclude stability in general. This means that trajectories can oscillate
around  $x^\star$ without ever converging. This  is intuitive since
$\eta(t)$ becomes very small for large $t$ and the system becomes almost
conservative---we refer to Ref. \cite{Franca:2018} for a more thorough stability
analysis of these systems.

Besides stability, it is possible to estimate
\emph{how fast} trajectories approach a minimum of $\F$.
This can also be done via a Lyapunov analysis, under certain convexity
assumptions on $\F$.   A function $\F$ is said to be \emph{convex} if
its gradient  $\nabla \F$
is maximal monotone, i.e., it obeys the inequality \eqref{mon}.
A function $\F$ is said to be
\emph{strongly convex} with parameter $\mu > 0$ if
it obeys a stronger condition:
$\langle \nabla \F(y) - \nabla\F(x) | y - x \rangle \ge \mu \| y - x\|^2$.
Let us mention some known rates of convergence
which follow from our results in \cite{Franca:2018b}.
For the gradient flow \eqref{gradflow} we have:
\begin{subequations}
\begin{align}
\F(x(t)) - \F^\star &= \bigO\big(t^{-1}\big) &\mbox{(convex),}
\label{gf_convex}\\
\| x(t) - x^\star\|^2 &= \bigO\big( e^{-\mu t}\big)
&\mbox{(strongly convex).} \label{gf_strongly_convex}
\end{align}
\end{subequations}
For the accelerated gradient flow \eqref{secode}
with constant damping,
$\eta(t) = \eta = \mbox{const.}$, we have:
\begin{subequations}
\begin{align}
\F(x(t)) - \F^\star &= \bigO\big(t^{-1}\big) &\mbox{(convex),}
\label{agf_const_convex}\\
\| x(t) - x^\star\|^2 &= \bigO\big( e^{-\sqrt{\mu} t}\big)
&\mbox{(strongly convex).} \label{agf_const_strongly_convex}
\end{align}
\end{subequations}
For the accelerated gradient flow \eqref{secode}
with decaying damping,
$\eta(t) = r/t$ where $r\ge 3$, we have:
\begin{subequations}
\begin{align}
\F(x(t)) - \F^\star &= \bigO\big(t^{-2}\big) &\mbox{(convex),}
\label{agf_decaying_convex} \\
\| x(t) - x^\star\|^2 &= \bigO\big( t^{-2 r / 3} \big)
&\mbox{(strongly convex).} \label{agf_decaying_strongly_convex}
\end{align}
\end{subequations}
We thus see that the accelerated gradient flow \eqref{secode} may converge
faster than the gradient flow \eqref{gradflow} in some situations.
For instance, Eq.~\eqref{agf_const_strongly_convex}  has
a $\sqrt{\cdot}$ improvement in the exponential
compared to \eqref{gf_strongly_convex},
while Eq.~\eqref{agf_decaying_convex} is an order of magnitude faster
compared to \eqref{gf_convex}.
Besides these rates, the stability of the system also plays a role.
However, one should note that these rates are upper bounds, and thus
may not always reflect the actual behavior of the system which may
be faster for a particular $\F$.

When $\F$ is quadratic we can solve Eqs.~\eqref{gradflow} and
\eqref{secode} exactly for some choices of damping $\eta(t)$.
Recall that this  captures
the behavior close to an isolated minimum:
$\F(x) - \F(x^\star) \approx
\tfrac{1}{2} (x-x^\star)^T  \nabla^2 \F(x^\star)  (x-x^\star) $.
We can change coordinates to a basis where the Hessian $\nabla^2 \F(x^\star)$
is diagonal so that the
components of the ODE become decoupled.
It is thus sufficient to consider the one-dimensional
case $\F(x) = \omega^2 x^2/2$.
Thus, the gradient flow \eqref{gradflow} has solution
\begin{equation} \label{sol_gf}
  x(t) = x_0 e^{-\omega^2 t},
\end{equation}
which agrees with the rate \eqref{gf_strongly_convex}.
The accelerated gradient flow \eqref{secode} with constant damping
has solution
\begin{equation} \label{sol_agf_const}
x(t) = x_0 e^{-\eta t / 2}
\left( \cosh\left( \xi t/2 \right) + (\eta/\xi) \sinh( \xi t / 2) \right) ,
\end{equation}
where $\xi \equiv \sqrt{\eta^2 - 4\omega^2}$ (we assumed
$\dot{x}_0 = 0$). This solution shows an exponential decay that
agrees with \eqref{agf_const_strongly_convex}.
Similarly, the accelerated gradient flow \eqref{secode} with a decaying
damping has solution
\begin{equation} \label{sol_agf_decaying}
x(t) = x_0 2^{\nu} \omega^{-\nu} \Gamma(\nu + 1)  t^{-\nu}
J_{ \nu }(\omega t), \qquad \nu \equiv \dfrac{r-1}{2},
\end{equation}
where $J_\nu$ is the Bessel function of the first kind (again,
$\dot{x}_0 = 0$).
%(we must
%remove $Y_\nu$, i.e. the Bessel function of the second kind,
%since it diverges at $t=0$).
A series expansion of $J_\nu$ for large $t$ reveals that
$J_\nu(t) \sim 1/\sqrt{t}$, which implies  $x(t) \sim t^{-r/2}$.
This is a power law, faster than the general upper bound
in \eqref{agf_decaying_strongly_convex}, however slower
than the exponential decay in \eqref{sol_agf_const}.

%%%%%%%%%%%%%%%%%%%%%%%%%%%
\begin{figure}
\centering
\includegraphics[scale=0.9]{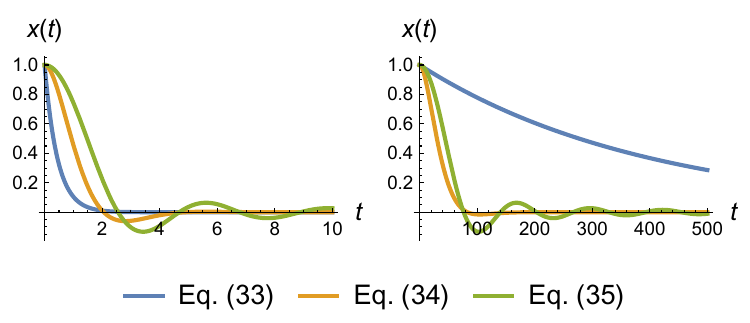}
\put(-200,100){\emph{(a)}}
\put(-55,100){\emph{(b)}}
\caption{Solutions of \eqref{gradflow} and
\eqref{secode} with $\F = (1/2) \omega^2 x^2$. For constant
damping \eqref{hbdamp}
we choose $\eta$ slightly below the critical value, and for decaying damping
\eqref{nagdamp} we choose $r = 3$.
\emph{(a)} $\omega = 1.5$.
\emph{(b)} $\omega=0.05$.
}
\label{solutions_fig}
\end{figure}
%%%%%%%%%%%%%%%%%%%%%%%%%%%

We illustrate the above solutions in Fig.~\ref{solutions_fig}.
Note that when $\omega > 1$ the gradient flow
tends to converge faster to equilibrium, however when $\omega < 1$
there
is a significant faster convergence of the accelerated gradient
flow in both cases.\footnote{Most challenging optimization
problems in machine learning tend to have $\omega  < 1$ since this
comes from a poor condition number of the objective function.}
In particular, the oscillations can
be better controlled with a constant damping.

The continuous-time dynamics is often easier to analyze
compared to the potentially complicated recurrence relations of a discrete-time
algorithm.
Thus, knowing the behavior of systems \eqref{gradflow} and \eqref{secode}
provide useful
insights to understand and design ``good'' optimization methods.
We expect that \emph{reasonable discretizations} of these systems
will reproduce their behavior, at least for sufficiently small choices of
the step size.
More precisely, it is a general result in numerical
analysis \cite{Hairer} that
if $\nabla \F$ is Lipschitz continuous with constant $L$,
a numerical integrator of order $p \ge 1$ has
\emph{global error}
$e_k \equiv  x(t_k) - x_k$ such that
$\| e_k \| \le C_t h^p$, where $t_k = k h$ and
$C_t = C (e^{L t} - 1)$ for some constant $C > 0$.
Thus, given a fixed simulation time $t$, one can control $e_k$
by making the step size $h$ sufficiently small. For our purposes,
the dissipative system presumably
converges to a minimum of $\F$ fast enough,
thus we expect that the simulation time will not be
excessively large which
helps controlling $C_t$---this
is a better scenario compared
to long-time simulations of conservative systems which are common, e.g., in
molecular dynamics and astrophysics.
Moreover, this estimate on $C_t$ is general and
pessimistic; a particular method may have a much smaller $C_t$.
In the following sections we will show that several well-known---but also
new---optimization methods are actually discretizations of the gradient flow
\eqref{gradflow}
or the accelerated gradient flow \eqref{secode}.

\section{Accelerated extensions of ADMM} \label{sec:admm}

\subsection{Balanced and rebalanced splitting}
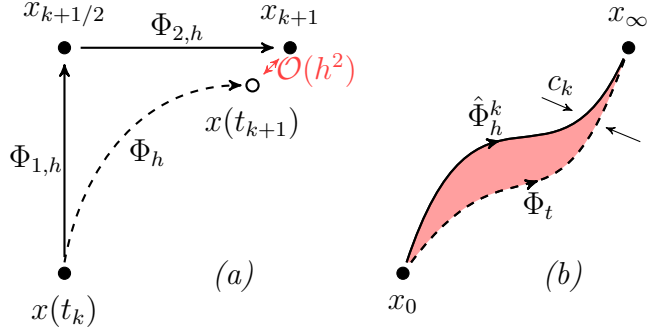
\begin{figure}
\centering
\tikzset{
    >=stealth',
    pil/.style={
        ->,
        thick,
        shorten <=2pt,
        shorten >=2pt,
    }
}
\tikzstyle directed=[postaction={decorate,decoration={markings,
    mark=at position .5 with {\arrow[scale=1.0]{>};}}}]

\begin{tikzpicture}[node distance=1cm, auto,]

\filldraw[black] (0,0) circle (2.3pt);
\filldraw[color=black,fill=white,thick] (2.5,2.5) circle (2.3pt);
\node[label=below:{$x(t_k)$}] (a) at (0,0) {};
\node[label=below:{$x(t_{k+1})$}] (b) at (2.5,2.5) {};
\path[->,thick,dashed] (a.north) edge[out=80,in=180] (b.west);

\filldraw[black] (0,3) circle (2.3pt);
\filldraw[black] (3,3) circle (2.3pt);
\node[label=above:{$x_{k+1/2}$}] (c) at (0,3) {};
\node[label=above:{$x_{k+1}$}] (d) at (3,3) {};
\path[pil] (a) edge (c);
\path[pil] (c) edge (d);

\path[<->,color=error1] (b) edge (d);

\node[color=error1] () at (3.35,2.7) {$\mathcal{O}(h^2)$};
\node[] () at (-0.38,1.5) {$\Phi_{1,h}$};
\node[] () at (1.5,3.25) {$\Phi_{2,h}$};
\node[] () at (1.1,1.6) {$\Phi_{h}$};

%%%%%%%%%%%%%%%%%

\filldraw[black] (4.5,0) circle (2.3pt);
\filldraw[black] (4.5+3,3) circle (2.3pt);
\node[label=below:{$x_0$}] (e) at (4.5,0) {};
\node[label=above:{$x_\infty$}] (f) at (4.5+3,3) {};

\draw[directed,thick,dashed,name path=A] (e) .. controls
                                            (6.0,2) and (6.5,0.2) .. (f);
\draw[directed,thick,name path=B] (e) .. controls
                                            (5.5,3.0) and (6.5,0.8) .. (f);
\tikzfillbetween[of=A and B]{error1, opacity=.5};
\draw[directed,thick,dashed] (e) .. controls
                                            (6.0,2) and (6.5,0.2) .. (f);
\draw[directed,thick] (e) .. controls
                                            (5.5,3.0) and (6.5,0.8) .. (f);

\node[] () at (6.3,0.9) {$\Phi_t$};
\node[] () at (5.6,2.1) {$\hat{\Phi}_{h}^{k}$};

\node[] (p11) at (6.9,2.1) {};
\node[] (p12) at (6.2,2.4) {};

\node[] (p21) at (7.05,2.0) {};
\node[] (p22) at (7.80,1.69) {};

\draw[<-] (p11) -- (p12);
\draw[<-] (p21) -- (p22);

\node[] () at (6.6,2.5) {$c_k$};

\end{tikzpicture}
\put(-170,20){\emph{(a)}}
\put(-45,20){\emph{(b)}}
\caption{\emph{(a)} Splitting approach
\eqref{composition}. In one time step the error between the
continuous and discrete trajectory is
$\| x(t_{k+1})  - x_{k+1}\| = \bigO(h^2)$.
\emph{(b)} The (re)balanced splitting introduces $c_k$ which
forces the discrete-time evolution to converge
to a critical point of the ODE: $\hat{\Phi}_h^k(x_0) =
\Phi_{t}(x_0) \to x_\infty$ as $k,t \to \infty$.}
\label{split_fig}
\end{figure}

Before making  connections with  ADMM we need to introduce
some ideas about splitting methods for ODEs.
Thus, consider the dynamical system
\begin{equation}
\label{first_order}
\dot{x} = A(x), \qquad A \equiv A_1 + A_2,
\end{equation}
where $A_1, A_2$ represent smooth and single-valued vector fields.
Suppose this is an intractable problem, i.e., the structure of $\F$ makes
the problem not amenable to a numerical procedure.
We denote the flow of \eqref{first_order}
by $\Phi_t$.
The idea  is to split the vector
field $A$ such that each individual system
\begin{equation}
\label{spliteq}
\dot{x} = A_1(x),
\qquad %\text{and} \qquad
\dot{x} = A_2(x) ,
\end{equation}
is integrable or has a feasible numerical
approximation. We denote their respective flows by
$\Phi_{1,t}$ and $\Phi_{2,t}$. For a step size $h > 0$,
it can be shown \cite{Hairer}
that the simplest composition
\begin{equation}
\label{composition}
\hat{\Phi}_h = \Phi_{2,h}\circ\Phi_{1,h}
\end{equation}
yields a first-order approximation, namely the \emph{local error} is
$\| \Phi_h(x) - \hat{\Phi}_h(x) \| = \bigO(h^2)$---see Fig.~\ref{split_fig}a
for an illustration.
However, in general, splittings such as \eqref{spliteq} \emph{do
not} preserve critical points of the original ODE.
The proposal of \cite{Speth:2013} is
to introduce a \emph{balance coefficient},
$c = c(t)$, and replace \eqref{spliteq} by
\begin{equation}
\label{balancedsplit}
\dot{x} = A_1(x) + c, \qquad  \dot{x} = A_2(x) - c.
\end{equation}
By appropriately choosing $c$ we can then  preserve critical points.
To see this, first suppose that $x_\infty$ is a critical point
of \eqref{first_order}, i.e.,
 $A_1(x_\infty) + A_2(x_\infty) = 0$. If $c_\infty$ obeys
\begin{equation} \label{cinfty}
c_\infty = \dfrac{1}{2}\big(A_2(x_\infty) - A_1(x_\infty)\big)
\end{equation}
then $x_\infty$ is also a critical point of both individual
Eqs. \eqref{balancedsplit}.
% since
%\begin{subequations}
%  \label{steadystate}
%\begin{align}
%\F_1(x_\infty) + c_\infty
%&= \tfrac{1}{2}\big(\F_1 + \F_2\big)(x_\infty) = 0, \\ % \qquad \text{and}  \\
%\F_2(x_\infty) - c_\infty
%&= \tfrac{1}{2}\big(\F_2 + \F_1\big)(x_\infty) = 0.
%\end{align}
%\end{subequations}
%where we have also used $\F_1(x_\infty) + \F_2(x_\infty) = 0$ from above.
Conversely, suppose  $x_\infty$ %$x_\infty = \lim_{t\to\infty}x(t)$
is a critical point of both individual Eqs.~\eqref{balancedsplit}. We  then
have
%, so that with $c_\infty \equiv \lim_{t\to\infty} c(t)$ it follows that
\begin{equation}
\begin{split} \label{cinfty.2}
c_\infty &= A_2(x_\infty) \\ &= - A_1(x_\infty) \\
&= \dfrac{1}{2}\big(A_2(x_\infty) - A_1(x_\infty)\big),
\end{split}
\end{equation}
implying that $x_\infty$ is also a critical point of
the original system \eqref{first_order}.
%From~\eqref{cinfty.2} and the fact that $x_\infty$ is stationary for both
%systems in~\eqref{balancedsplit} gives
%\begin{subequations}\label{steadystate.all}
%  \label{steadystate.2}
%\begin{align}
%0 &= \F_1(x_\infty) + c_\infty
%   = \tfrac{1}{2}\big(\F_1 + \F_2\big)(x_\infty), \\ % \qquad \text{and}  \\
%0 &= \F_2(x_\infty) - c_\infty
%   = \tfrac{1}{2}\big(\F_2 + \F_1\big)(x_\infty),
%\end{align}
%\end{subequations}
%so that both equations in~\eqref{steadystate.all} imply that $x_\infty$ is
% stationary for~\eqref{first_order}.
%This shows that if a point is simultaneously stationary for both
% dynamical systems in~\eqref{balancedsplit}, then it is also stationary
% for~\eqref{first_order}. %, i.e., that stationary states are preserved.
This can be implemented with
$c_{k+1} = \tfrac{1}{2}\big(A_2(x_k) - A_1(x_k) \big)$,
together with suitable discretizations of
Eqs.~\eqref{balancedsplit}---see Fig.~\ref{split_fig}b for an illustration.
However, this approach requires explicit
computation of the vector fields $A_i$'s.  In optimization this means
computing gradients, which may not be available.
To address this issue we consider a related approach.

%%%%%%%%%%%%%%%%%%%%%%
%%%%%%%%%%%%%%%%%%%%%%%%%%%%%%%%%%%%

The \emph{rebalanced splitting}  \cite{Speth:2013} is particularly
well-suited for the \emph{implicit} discretizations we have in mind.
We thus integrate
$\dot{x} = A_1(x) + c_k$
over the interval $[t_k, t_k+h]$,
with initial condition
$x(t_k) = x_k$,
 to obtain the intermediate point
$x_{k+1/2}$. Then we
integrate
$\dot{x} = A_2(x) - c_k$ over the same interval,
with initial condition $x(t_k) = x_{k+1/2}$, to obtain
the endpoint $x_{k+1}$. Note that $c_k$ is kept fixed during this procedure.
Thus,
\begin{subequations}
\label{rebal}
\begin{align}
x_{k+1/2} &= x_k +
\int_{t_k}^{t_k + h}\big(A_1(x(t))+c_k\big)dt,  %\qquad \text{and}
\label{rebal1}\\
x_{k+1} &= x_{k+1/2} +
\int_{t_k}^{t_k + h}\big(A_2(x(t))-c_k\big)dt. \label{rebal2}
\end{align}
\end{subequations}
In light of \eqref{cinfty.2}, two reasonable ways of
computing $c_{k}$ are  %namely as
given by the \emph{average} of either $\tfrac{1}{2}(A_2 - A_1)$ or $A_2$.
We choose the latter---as we will see this allow us to derive
ADMM---which with \eqref{rebal} yields
\begin{equation}
\label{crebal}
\begin{split}
%c_{k+1}
%&= \dfrac{1}{h}\int_{t_k}^{t_k + h}\!\! \dfrac{ \F_2(x(t)) - \F_1(x(t)) }{2} dt \nonumber
%\\ &= c_k + \dfrac{1}{h}\left( \dfrac{x_{k+1} + x_k}{2} - x_{k+1/2} \right) ,  \label{crebal.1} \\
c_{k+1}
&= \dfrac{1}{h}\int_{t_k}^{t_k+h} A_2(x(t)) dt \\
&= c_k + h^{-1}\left(x_{k+1}-x_{k+1/2} \right).
\end{split}
\end{equation}
In contrast to the previous case, we now
need not compute $A_i$ explicitly.
Our derivation above is slightly different than the one
in Ref. \cite{Speth:2013}.

\subsection{Deriving extensions of ADMM}
We are now  in a position to show how generalizations of ADMM
emerge from
such an approach. We focus on problem
\begin{equation}\label{optimization2}
\min_{x \in \mathbb{R}^n} \F(x), \qquad \F=\F_1+\F_2+\F_3 ,
\end{equation}
and moreover on discretizations of the accelerated gradient flow
\eqref{secode} since
discretizations of the
gradient flow \eqref{gradflow} can be recovered
as particular cases; recall the discussion in Sec.~\ref{sec:high_friction}.

With a balance coefficient $c=c(t)$ we write
\eqref{secode} as
\begin{subequations}
\label{secbal}
\begin{align}
\dot{x} &= p , \\
\dot{p} &= \underbrace{-\eta(t) p - \nabla \F_1(x) -\nabla \F_3(x)}_{A_1} +c
\underbrace{-\nabla \F_2(x)}_{A_2} - c  .
\end{align}
\end{subequations}
Splitting as indicated, and further combining
%\eqref{secbal} yields
%\begin{equation}
%\dot{x} = v, \qquad  \dot{v} = -\eta(t) v - \nabla \F_1(x) - \nabla \F_3(x) + c ,
%\end{equation}
the resulting  equations, we obtain two independent ODEs:
%\begin{equation}
%\dot{x} = v \qquad  \dot{v}  = -\nabla \F_2(x) - c .
%\end{equation}
%It is convenient to combine these individual systems back into second order form,
\begin{align}
\ddot{x} +  \eta(t) \dot{x}  &= -\nabla \F_1(x) - \nabla \F_3(x) + c ,
\label{admm_split1} \\
\ddot{x}  &= -\nabla \F_2(x) - c . \label{admm_split2}
\end{align}
%so as to discretize them using similar approach as in subsection~\ref{sec:building_blocks}.
%Note that we will replace the step size $h^2 \to h$
%as already done in \eqref{a_proximal_point}.
Using \eqref{eulersec}, and redefining the step size $h^2 \mapsto h$ as
discussed in \eqref{a_proximal_point},
a semi-implicit discretization of \eqref{admm_split1} is
\begin{equation} \label{admm_disc1}
\dfrac{x_{k+1/2} - \hat{x}_k}{h} = - \nabla \F_1(x_{k+1/2}) -
\nabla \F_3(\hat{x}_k) + c_k .
\end{equation}
This can  be solved
with the resolvent \eqref{resolvent} to obtain
\begin{equation} \label{admm_up1}
x_{k+1/2} = J_{h \partial \F_1}\big( \hat{x}_k -
h \nabla \F_3(\hat{x}_k) + h c_k \big).
\end{equation}
We now discretize \eqref{admm_split2} as
\begin{equation} \label{xtilde}
\dfrac{\tilde{x}_{k+1} -2 x_{k+1/2} + \hat{x}_k}{h} =
- \nabla \F_2(x_{k+1}) - c_k
\end{equation}
where
\begin{equation}\label{xtilde1}
\tilde{x}_{k+1} \equiv x_{k+1} + (x_{k+1/2} - \hat{x}_k).
\end{equation}
Note that $\tilde{x}_{k+1}$ is related to  $x_{k+1}$ via the
``momentum'' term $(x_{k+1/2} - \hat{x}_k)$ based on the
first splitting.\footnote{$\tilde{x}_{k+1}$ is
slightly further away from
$x_{k+1}$ which makes the algorithm ``look ahead'' and implicitly introduces
dependency on the curvature of $\F_2$. Although the introduction
of $\tilde{x}_{k+1}$ may seem artificial, it will be justified below when
we compute the error in approximating the continuous trajectory.}
%%%%%%%
%\begin{equation} \label{admm_disc2}
%$x_{k+1} - x_{k+1/2} = -h \nabla \F_2(x_{k+1}) - h c_k$
%\end{equation}
With  this and the resolvent we obtain
\begin{equation} \label{admm_up2}
  x_{k+1} = J_{h \partial \F_2}\big( x_{k+1/2} - h c_k \big).
\end{equation}
The balance coefficient follows readily from \eqref{crebal}:%%%%%%
%%%%%%%%%%%%%%%%%
\footnote{To justify that $h^2\mapsto h$ does not change this, note that
with $A_2  = -\nabla \F_2$ in \eqref{crebal} an implicit discretization
corresponds to approximating the integral by its upper limit, thus
$c_{k+1} = \tfrac{1}{h}\int_{t_k}^{t_k+h} B dt  \approx
- \nabla \F_2(x_{k+1})$.
Using \eqref{xtilde}--\eqref{xtilde1} yields \eqref{admm_up3}.}
%%%%%%%%%%%%%%%%%
\begin{equation}\label{admm_up3}
c_{k+1} = c_k + h^{-1}\left( x_{k+1}-x_{k+1/2} \right).
\end{equation}
Combining the above steps results into
a family of accelerated extensions of ADMM  summarized in
Algorithm~\ref{agenadmmrebal}.

\begin{algorithm}
\DontPrintSemicolon

    Choose step size $h$ and damping function $\gamma_k$\;
    Initialize $x_0$, $\hat{x}_0$ and $c_0=0$\;

    \For{$k=0,1,\dotsc$}{
        \leftskip 15pt

        $x_{k+1/2} = \prox_{h \F_1}\left(\hat{x}_k -
                     h \nabla \F_3(\hat{x}_k) + h c_k\right)$\;
        $x_{k+1}  = \prox_{h \F_2}( x_{k+1/2} - h c_k)$\;
        $c_{k+1} = c_k + h^{-1}\left( x_{k+1}-x_{k+1/2}\right)$\;
        $\hat{x}_{k+1} = x_{k+1} + \gamma_{k+1}(x_{k+1}-x_k)$\;

    }

\caption{Family of accelerated extensions of ADMM for problem
\eqref{optimization2}.\label{agenadmmrebal}}
\end{algorithm}

Let us stress some important aspects of Algorithm~\ref{agenadmmrebal}.
The standard ADMM \cite{Gabay:1976,Glowinsky:1975,Boyd:2011}
corresponds to the particular case where
$\F_3=0$ and no acceleration is used, i.e., $\gamma_k=0$.
Thus, Algorithm~\ref{agenadmmrebal} not only generalizes ADMM to handle
problems in the form \eqref{optimization2}
but also includes acceleration with arbitrary damping functions $\eta(t)$.
The so-called dual vector in ADMM,  originally obtained as a
Lagrange multiplier \cite{Boyd:2011}, is here represented by the balance
coefficient $c_k$ and thus acquires a new meaning:
its role is to preserve critical points of the ODE.

When decaying damping \eqref{nagdamp} is chosen and $\F_3=0$,
Algorithm \ref{agenadmmrebal} is similar to the so-called ``fast ADMM''
\cite{Goldstein:2014}. They differ in that the latter also  ``accelerates''
the dual variable $c_k$. Connections between fast ADMM and ODEs
was considered recently by us in Refs. \cite{Franca:2018,Franca:2018b} and
also corresponds to Eq. \eqref{secode}, however  the discretization
is not a rebalanced splitting within the above framework.

In light of the discussion in Sec.~\ref{sec:resolvent}, it is clear
that Algorithm~\ref{agenadmmrebal} can be generalized to account for monotone
operators by the replacements
$\prox_{h \F_1} \to J_{h A_1}$ and $\prox_{h \F_2} \to J_{h A_2}$.

Finally,  although we focused on the accelerated gradient flow
\eqref{secode}, analogous procedure
applies to the gradient flow \eqref{gradflow} and leads to
Algorithm~\ref{agenadmmrebal} with $\gamma_k=0$.%%%%%%%%%%
%%%%%%%%%%%%%%%%%%%%%%%%%
\footnote{We have $\dot{x} = - \nabla (\F_1 + \F_3)(x) + c$ and
$\dot{x} = - \nabla \F_2(x) - c$. For the former,
$x_{k+1/2} - x_k = - h \nabla\F_1(x_{k+1/2}) - h\nabla \F_3(x_k) + h c_k$,
whereas
for latter, $x_{k+1} - x_{k+1/2} = - h\nabla \F_2(x_{k+1}) - h c_{k}$.
Using the resolvent and
\eqref{crebal} yield Algorithm~\ref{agenadmmrebal} with $\gamma_k=0$.}
%%%%%%%%%%%%%%%%%%%%%%%%

\subsection{Order of accuracy} \label{sec:accuracy_aadmm}
Next, we show that the above discretization is justified since it is a
first-order approximation to the continuous trajectory,
i.e., $\| \Phi_h(x) - \hat{\Phi}_h(x)\| = \bigO(h^2)$.
From the definition of the resolvent \eqref{resolvent} we have that
$y = J_{ h \nabla \F}(x)$ if and only if
$y  = x - h \nabla \F(y)$, thus
\begin{equation} \label{proxappro}
y %J_{ h \nabla \F}(x) \\
=  x - h \nabla \F(x - h \nabla \F(y))
= x - h \nabla \F(x) + \Order(h^2) .
\end{equation}
%where the
%$\|\nabla f(y)\|$
%that normally appears in the $\Order$ term is  suppressed because it is bounded independently of $\lambda$ for all $\lambda$ on a compact set as a consequence of $y = J_{\lambda \nabla f}(x)$ and Assumption~\ref{assump2}. (A similar convention is used in~\eqref{limm} and~\eqref{limmm} below.)
This relation implies the following approximations, valid up to $\bigO(h^2)$,
for the updates in Algorithm~\ref{agenadmmrebal}:
\begin{align}
x_{k+1/2} &\approx \hat{x}_k - h \nabla \F_3(\hat{x}_k) + h c_k -
h \nabla \F_1(\hat{x}_k) , \label{limm} \\
x_{k+1} &\approx
%x_{k+1/2} - h c_k - h \nabla \F_2(x_{k+1/2}) + \Order(h^2) \\
\hat{x}_k - h \nabla \F(\hat{x}_k) .
\label{limmm}
\end{align}
%where to derive the second equality we used~\eqref{limm} and a Taylor expansion of $\nabla g$.
Recall \eqref{xhat}, namely $\gamma_k = 1- \sqrt{h} \eta(t_k) $
since we redefined $h^2\mapsto h$. Thus,
% for constant damping $(\eta(t)=r)$, while
%\begin{equation}
%\gamma_k = \dfrac{k}{k+r} = 1 - \dfrac{r}{k+r} = 1 - \dfrac{r h }{t_k} \left(1+ \dfrac{rh}{t_k}\right)^{-1}
%= 1 -  \eta(t_k) h + \Order(h^2)
%\end{equation}
%for decaying damping ($\eta(t) = r/t$).
%Thus, in either case, we conclude that
\begin{equation}
\label{xhat_ord}
%\begin{split}
\hat{x}_k = x_k + \big( 1 - \eta(t_k) \sqrt{h} \big) \sqrt{h} p_k   %+ \Order(h^{3/2}) %\\
%&
%&= x_k + \Order(\sqrt{h})
%\end{split}
\end{equation}
where
\begin{equation}\label{velocity}
p_k \equiv \dfrac{ x_k - x_{k-1} }{ \sqrt{h} } .
\end{equation}
From \eqref{xhat_ord} and \eqref{limmm}, and now restoring the original
step size ($h \mapsto h^2$), we conclude that
%to $\bigO(h^2)$ it holds that
\begin{align}
p_{k+1} &\approx p_k - h \eta(t_k) p_k - h \nabla \F(x_k),
\label{admm_first_final1} \\
x_{k+1} &= x_k + h p_{k+1} \approx x_k + h p_k .
\label{admm_first_final2}
\end{align}
Finally,  the evolution of Eq. \eqref{secode} in one time step gives
\begin{align}
p(t+h) &= p(t) + h \dot{p}(t) + \Order(h^2)  \nonumber \\
&\approx p(t) - h \eta(t) p(t) - h \nabla \F(x(t)) , \\
x(t+h) &= x(t) + h \dot{x}(t) + \Order(h^2) \nonumber \\
&\approx x(t) + h p(t) .
\end{align}
Comparing with \eqref{admm_first_final1}--\eqref{admm_first_final2}
implies that the algorithm's state simulates the continuous trajectory
up to an error $\bigO(h^2)$, therefore
%\begin{equation} \label{admm_first_order_final}
%v(t_{k+1}) =   v_{k+1} + \Order(h^2), \quad %\mbox{and} \qquad
%x(t_{k+1}) = x_{k+1} + \Order(h^2) ,
%\end{equation}
the discretization is first-order accurate.
Similar conclusion holds for the nonaccelerated algorithm in relation
to the gradient flow \eqref{gradflow}.

We mention a subtlety when $\F$ is nonsmooth, or when
considering monotone operators more generally.
A  crucial step in the above derivation was the Taylor
approximation of the resolvent \eqref{proxappro}.
For a maximal monotone operator $A$,
in the most general case only a slightly weaker approximation
is available \cite[Remark 23.47]{BauschkeBook}:
\begin{equation} \label{taylor_res}
J_{h A}  = I - h A_0  + \smallO(h)
\end{equation}
where $A_0 x = \lim_{\mu \downarrow 0}A_\mu x$
(see Sec.~\ref{sec:resolvent}).
The previous arguments still hold true, however due
to \eqref{taylor_res}, and assuming that we can expand
$A_0(x + \Order(h)) = A_0(x) + \Order(h)$, the local error is now
$\smallO(h)$ instead of $\Order(h^2)$.
It is important to note that this is a consequence of the nonsmoothness
of $\F$,  or the multivaluedness of $A$,
and not of the discretization procedure.
%In other words,
%within this approach there is nothing one can do about it.
This comment applies to all cases considered in this paper.

\section{Accelerated extensions of Davis-Yin} \label{sec:davis_yin}

\subsection{Discretization}
We now introduce accelerated extensions of Davis-Yin \cite{Davis:2017}.
This time we split the system \eqref{secode} without
a balance coefficient, namely we choose vector fields
\begin{align}
A_1(x) &= -\eta(t) \dot{x} - \nabla \F_1(x), \\ %\text{and}  \qquad
A_2(x) &= - \nabla \F_2(x) - \nabla \F_3(x).
\end{align}
Instead of Eqs.~\eqref{admm_split1} and \eqref{admm_split2} we now
obtain
\begin{align}
\ddot{x} + \eta(t)\dot{x} &= -\nabla \F_1(x),
\label{ady_flow1} \\
\ddot{x} &= -\nabla \F_2(x) - \nabla \F_3(x).
\label{ady_flow2}
\end{align}
Using \eqref{eulersec} (and redefining $h^2 \mapsto h$),
an implicit discretization of Eq. \eqref{ady_flow1} is
\begin{equation} \label{ady.disc1}
\dfrac{x_{k+1/4} - \hat{x}_k}{h} = -\nabla \F_1(x_{k+1/4}),
\end{equation}
which with the resolvent \eqref{resolvent}  gives
\begin{equation} \label{ady.u1}
x_{k+1/4} \equiv \Phi_{1,h}(\hat{x}_k) =  J_{h \partial \F_1}(\hat{x}_k) .
\end{equation}
Next, to ``inject momentum'' in the direction of $\nabla \F_1$,
we define the ``translation operator''
\begin{equation} \label{translation}
\mathcal{T}_h(z) \equiv z - h \nabla \F_1(x_{k+1/4}) .
\end{equation}
The next point is thus obtained as
\begin{equation} \label{ady.u2}
x_{k+1/2} \equiv \mathcal{T}_{h}(x_{k+1/4}) =  2x_{k+1/4} - \hat{x}_k.
\end{equation}
%This last equation
%is a gradient descent step on $f$.
A semi-implicit discretization of Eq. \eqref{ady_flow2} is
\begin{equation}\label{adyfirst}
\dfrac{x_{k+3/4} - 2x_{k+1/4} + \hat{x}_k}{h}
= - \nabla \F_2(x_{k+3/4}) - \nabla \F_3(x_{k+1/4}) ,
\end{equation}
which can again be solved with the resolvent \eqref{resolvent} to
obtain
\begin{equation} \label{ady.u3}
\begin{split}
x_{k+3/4} &\equiv \Phi_{2,h}(\hat{x}_k)  \\
&= J_{h \partial \F_2}\left(x_{k+1/2} - h \nabla \F_3(x_{k+1/4})
\right).
\end{split}
\end{equation}
%Let us make some observations about the operator $\mathcal{T}_{h}$ in \eqref{ady.u2}.
%First, it translates the vector $x_{k+1/4}$ twice the distance away from $\hat{x}_k$, i.e., $x_{k+1/2}$ satisfies
%$x_{k+1/2} - \hat{x}_k = 2 (x_{k+1/4} - \hat{x}_k)$.
%Second, it produces the same point as the Cayley operator acting on $\hat{x}_k$ since
%$\mathcal{T}_h(x_{k+1/4}) = (2 J_{\lambda \nabla f} - I)(\hat{x}_k) \equiv C_{\lambda \nabla f}(\hat{x}_k)$ (see Lemma~\ref{caleylemma}).  To see that $\mathcal{T}_h$ is a well-defined function of $x_{k+1/4}$, note that we can also write
%\begin{equation}
%\mathcal{T}_h(x_{k+1/4})= \big(2I - J_{\lambda \nabla f}^{-1}\big)(x_{k+1/4}) =
%(I - \lambda \nabla f)(x_{k+1/4}),
%\end{equation}
%where we used the definition of the resolvent \eqref{resolvent} and \eqref{ady.u1}.
%Thus, the action of $\mathcal{T}_{h} = I - \lambda \nabla f$ on $x_{k+1/4}$ is
%a gradient descent step, and has an inverse  $\mathcal{T}_h^{-1}$ (whose specific form is not necessary for our purposes).
%The introduction of \eqref{ady.u2} may seem artificial, although completely allowed, and we will comment more about this shortly.
Finally,  applying the inverse
$\mathcal{T}^{-1}_h(z) \equiv z + h \nabla \F_1(x_{k+1/4} )$
and using \eqref{ady.disc1} we obtain
%for any vector $z$, hence
%we just need to discount the translation by subtracting $x_{k+1/4}-\hat{x}_k$ (this is much easier than trying to find $\mathcal{T}_h^{-1}$ explicitly), thus
\begin{equation} \label{ady.u4}
x_{k+1} \equiv \mathcal{T}^{-1}_{h}(x_{k+3/4})
=  x_{k+3/4}-(x_{k+1/4}-\hat{x}_k).
\end{equation}
We collect the above steps into Algorithm~\ref{ady}.
The original
Davis-Yin method \cite{Davis:2017} is
recovered by setting $\gamma_k = 0$. Such a case corresponds
to an overdamped limit---recall the discussion
of Sec.~\ref{sec:high_friction}---which is indeed a
discretization of the gradient
flow \eqref{gradflow}, as can also be easily verified by repeating the above
procedure to this simpler case.

\begin{algorithm}
\DontPrintSemicolon

    Choose step size $h$ and damping function $\gamma_k$\;
    Initialize $x_0$ and $\hat{x}_0$\;

    \For{$k=0,1,\dotsc$}{
        \leftskip 15pt

        $x_{k+1/4} = \prox_{h \F_1}(\hat{x}_k)$\;
        $x_{k+1/2} =  2x_{k+1/4} - \hat{x}_k$\;
        $x_{k+3/4} = \prox_{h \F_2}\big(x_{k+1/2} -
                    h \nabla \F_3(x_{k+1/4})\big)$\;
        $x_{k+1} = \hat{x}_k + x_{k+3/4} - x_{k+1/4}$\;
        $\hat{x}_{k+1} = x_{k+1} + \gamma_{k+1}(x_{k+1}-x_k)$\;

    }

\caption{Family of accelerated extensions of Davis-Yin (DY) for
problem \eqref{optimization2}.\label{ady}}
\end{algorithm}

Algorithm \ref{ady} is equivalent to the fixed point iteration
$x_{k+1}=\hat{\Phi}_h(\hat{x}_k)$ with
\begin{equation}
\label{dycompos}
\hat{\Phi}_h \equiv
\mathcal{T}^{-1}_h\circ\Phi_{2, h}\circ\mathcal{T}_h\circ\Phi_{1, h} ,
\end{equation}
where these individual  maps are defined in
Eqs. \eqref{ady.u1}, \eqref{translation},
and \eqref{ady.u3}.
Thus, the translation operator  $\mathcal{T}_h$ is
actually a ``preprocessor map,'' which is a common technique
in numerical analysis  \cite{Hairer}.
The discretization associated to Davis-Yin can be summarized
diagrammatically:
\begin{equation}
  \begin{tikzcd}[column sep=large, row sep=huge]
    \hat{x}_k \arrow[thick]{r}{\Phi_{1,h}} \arrow[dashed,thick]{rd} & x_{k+1/4}
    \arrow[thick]{r}{\mathcal{T}_h} & x_{k+1/2} \arrow[thick]{d}{\Phi_{2,h}} \\
     & x_{k+1} \arrow[leftarrow,thick]{r}{\mathcal{T}_{h}^{-1}} & x_{k+3/4}
  \end{tikzcd}
\end{equation}

\subsection{Order of accuracy} \label{sec:accuracy_ady}

Using the expansion \eqref{proxappro} we can approximate the updates of
Algorithm~\ref{ady} up to $\Order(h^2)$:
\begin{subequations}
\begin{align}
x_{k+1/4} &\approx \hat{x}_k - h  \nabla \F_1(\hat{x}_k), \\
x_{k+1/2} &\approx \hat{x}_k - 2 h \nabla \F_1(\hat{x}_k), \\
x_{k+3/4} &\approx \hat{x}_k - 2 h \nabla \F(\hat{x}_k), \\
x_{k+1}  &\approx \hat{x}_k  - h \nabla \F(\hat{x}_k). \label{last_exp}
\end{align}
\end{subequations}
But \eqref{last_exp} is exactly the same as \eqref{limmm},
thus the remaining steps of the argument follow as before,
implying that Algorithm~\ref{ady} is also a first-order integrator
to the accelerated gradient flow \eqref{secode}. The same holds true
for standard Davis-Yin
($\gamma_k=0$) in relation to the gradient flow~\eqref{gradflow}.

\subsection{Preserving critical points}

Since Algorithm~\ref{ady} arises from a splitting that is not balanced,
it is not a priori obvious if critical points of the underlying ODE are
preserved. We now show that this is indeed the case.
We can write the operator \eqref{dycompos} as
\begin{equation}
\label{dycaley}
%\begin{split}
\hat{\Phi}_h =
I + J_{h \partial \F_2}\circ\left(2 J_{h \partial \F_1} - I - h \nabla \F_3
\circ J_{ h \partial \F_1}\right)  - J_{ h \partial \F_1}.
%
%&= I + J_{ h \partial \F_2 } \circ\left(  C_{ h \partial \F_1 } - h \nabla \F_3 \circ
%J_{ h \partial \F_1 } \right) - \dfrac{1}{2}(C_{ h \partial \F_1 } + I) \\
%&= \dfrac{1}{2}I - \dfrac{1}{2} C_{ h \partial \F_1 } +  \dfrac{1}{2}\left(C_{ h \partial \F_2 } + I\right)\circ
%\left( C_{ h \partial \F_1 }  - h \nabla \F_3 \circ J_{ h \partial \F_1 } \right) \\
%&= \dfrac{1}{2}I +\dfrac{1}{2}C_{ h \partial \F_2 }\circ\left( C_{h \partial \F_1 }
%- h  \nabla \F_3  \circ J_{h  \partial \F_1 } \right) - \dfrac{1}{2} h  \nabla \F_3 \circ
%J_{h \partial \F_1} .
%\end{split}
\end{equation}
%\end{widetext}
Assuming the algorithm converges, we must have a fixed point equation
$x_{\infty} = \hat{\Phi}_h(x_{\infty})$.  We thus need to show that this
equation generates critical points of the ODE. To this end,
let  $x^\star$ be such a critical point, i.e.,
\begin{equation} \label{crit}
(\nabla \F_1 + \nabla \F_2 + \nabla \F_3)(x^\star) = 0.
\end{equation}
This is equivalent to
$(I + h \nabla \F_2)(x^\star) = (I - h \nabla \F_1 -h\nabla \F_3)(x^\star)$,
and with the aid of the resolvent \eqref{resolvent} it can be written as
\begin{equation}
  x^\star = J_{h \partial \F_2} \circ
  (I - h \nabla \F_1 - h \nabla \F_3)(x^\star).
\end{equation}
Using the identity
\begin{equation}
(2J_{h \partial \F} - I)  \circ (I + h \nabla \F) = I - h \nabla \F
\end{equation}
we thus have
\begin{equation}
x^\star = J_{h \partial \F_2} \circ
\left[ (2J_{h \partial \F_1}-I)\circ(I+h \nabla \F_1) -
h \nabla \F_3 \right](x^\star).
\end{equation}
Define $x_\infty \equiv (I + h \nabla \F_1)(x^\star)$, i.e.,
$x^\star = J_{h \partial \F_1}(x_\infty)$.
The above equation then yields
\begin{equation}
J_{h \partial \F_1}(x_\infty) = J_{h \partial \F_2} \circ
    \left[ 2J_{h \partial \F_1}-I - h \nabla \F_3 \circ
    J_{h \partial \F_1} \right](x_\infty)
\end{equation}
which is equivalent to $x_\infty = \hat{\Phi}_h(x_\infty)$
according to \eqref{dycaley}.
Therefore, critical points of \eqref{crit} yield fixed points of the operator
\eqref{dycaley} and vice-versa.  This shows that
Algorithm~\ref{ady} preserves critical points of the underlying ODE.

\subsection{Accelerated extensions of Douglas-Rachford}
\label{sec:douglas_rachford}

Douglas-Rachford (DR) \cite{Douglas:1956,Lions:1979} is
recovered from Algorithm~\ref{ady}
in the particular case where $\gamma_k=0$ and
$\F_3 = 0$.
Therefore, in the case where $\F_3 = 0$ but $\gamma_k \ne 0$,
Algorithm~\ref{ady} contains accelerated extensions of
Douglas-Rachford---the case with decaying damping was considered
in \cite{Patrinos:2014}.
From the previous arguments, we know that such methods are
discretizations of the accelerated gradient flow \eqref{secode},
whereas the standard Douglas-Rachford is
a discretization of the gradient flow \eqref{gradflow}.
Moreover, such discretizations preserve critical points and are
first-order integrators.

\subsection{Accelerated extensions of forward-backward}
\label{sec:forward_backward}

The forward-backward method \cite{Lions:1979} is recovered from
Algorithm~\ref{ady} when $\gamma_k=0$ and $\F_1 = 0$.
Thus, when $\F_1=0$ but $\gamma_k \ne 0$, Algorithm~\ref{ady} reduces
to
\begin{equation} \label{fbmethod}
\begin{split}
x_{k+1} &= \prox_{h  \F_2 }\big( \hat{x}_k -
h \nabla \F_3(\hat{x}_k) \big), \\
\hat{x}_{k+1} &= x_{k+1} + \gamma_{k+1} (x_{k+1} - x_k).
\end{split}
\end{equation}
The decaying damping case \eqref{nagdamp} was
considered in \cite{Beck:2009}.
From an ODE perspective, the above discretization
is not a splitting method but rather a direct
semi-implicit discretization of Eq. \eqref{secode}.
Anyhow, our previous arguments show that such accelerated variants of
forward-backward are first-order integrators of this ODE
and preserve critical points---the same holds true for the
standard forward-backward ($\gamma_k=0$) in relation to the
gradient flow~\eqref{gradflow}.

\section{Accelerated extensions of Tseng's splitting} \label{sec:tseng}

The last proximal method remaining to be considered is the
forward-backward-forward or Tseng's splitting \cite{Tseng:2000}.
%, which
%is a  ``slight perturbation'' of
%forward-backward.
%We anticipate that in practice, we did not see improvements of (accelerated) Tseng's versus (accelerated)
%forward-backward, however we still consider this case for
%completeness and theoretical reasons.
Thus, consider  Eq. \eqref{secode} with $\F_1=0$ and written as
%\begin{subequations}
\begin{align}
%\dot{x} &= p, \\
\dot{p} &=  \underbrace{-\eta(t) p -\nabla \F_2(x) -
\nabla \F_3(x)}_{A_1} + \underbrace{\nabla \F_3(x) - \nabla \F_3(x)}_{A_2}.
\end{align}
%\end{subequations}
Note that $A_2 = 0$, however in a discretization
there might be numerical inaccuracies which
introduces a kind of ``perturbation'' on top of the forward-backward method,
which arises from the first component of this system.
Splitting as indicated yields
\begin{align}
  \ddot{x} + \eta(t) \dot{x} &= -\nabla \F_2(x)-\nabla \F_3(x),
  \label{tseng_split1} \\
  \ddot{x} &= \nabla \F_3(x)-\nabla \F_3(x). \label{tseng_split2}
\end{align}
A semi-implicit
discretization of the first equation yields
\begin{equation} \label{tseng.u1}
x_{k+1/2} =
J_{h \partial \F_2 }\left( \hat{x}_k - h \nabla \F_3 (\hat{x}_k) \right),
\end{equation}
which is the forward-backward method \eqref{fbmethod}.
Eq. \eqref{tseng_split2} can be discretized as
$\tilde{x}_{k+1} - 2 x_{k+1/2} +
\hat{x}_k = h \nabla \F_3(\hat{x}_k) - h \nabla \F_3(x_{k+1/2})$, where
$\tilde{x}_{k+1}$  is given by  \eqref{xtilde1}.
%as already introduced for ADMM.
Thus,
\begin{equation}\label{tseng.u2}
x_{k+1} = x_{k+1/2} -
h \big( \nabla \F_3 (x_{k+1/2}) - \nabla \F_3(\hat{x}_k) \big).
\end{equation}
Therefore, we derived Algorithm~\ref{atseng}.

\begin{algorithm}
\DontPrintSemicolon

    Choose step size $h$ and damping function $\gamma_k$\;
    Initialize $x_0$ and $\hat{x}_0$\;

    \For{$k=0,1,\dotsc$}{
        \leftskip 15pt

        $x_{k+1/2} = \prox_{h  \F_2}\left(\hat{x}_k- h
                        \nabla \F_3(\hat{x}_k)\right)$\;
        $x_{k+1}  = x_{k+1/2} - h \left( \nabla \F_3(x_{k+1/2})-
                            \nabla \F_3(\hat{x}_k) \right)$\;
        $\hat{x}_{k+1} = x_{k+1} + \gamma_{k+1}(x_{k+1}-x_k)$\;

    }

\caption{Family of accelerated extensions of Tseng's method.\label{atseng}}
\end{algorithm}

The original Tseng's splitting is recovered with $\gamma_k=0$, in which case
it is a discretization of the gradient flow \eqref{gradflow}.
Due to Eq. \eqref{tseng_split2} we expect to have
a ``contraction'' on the acceleration which indicates that
Algorithm~\ref{atseng} may
be slower than \eqref{fbmethod} (this was actually observed
in our experiments and this method tends to be slower than forward-backward).

In a similar way as already done in
Secs.~\ref{sec:admm} and \ref{sec:davis_yin},
through Taylor expansions it is straightforward to show that the above
discretization is first-order accurate.
%---we omit this derivation for
%conciseness.

We can also show that the above discretization preserves critical points.
Indeed, Algorithm~\ref{atseng} is equivalent to iterations
$x_{k+1} = \hat{\Phi}_h(\hat{x}_k)$ with
\begin{equation} \label{tsengop}
\hat{\Phi}_h
= (I - h \nabla \F_3)\circ
  J_{ h \partial \F_2 }\circ(I - h  \nabla \F_3) + h \nabla \F_3.
\end{equation}
Assuming the algorithm converges,
$x_\infty = \hat{\Phi}_h(x_\infty)$, i.e.,
\begin{equation}
(I-h\nabla\F_3)(x_\infty) =
(I-h\nabla\F_3)\circ J_{h\partial\F_2}\circ(I-h\nabla\F_3)(x_\infty).
\end{equation}
Moreover, assuming that $h$ is sufficiently small so that the inverse
$(I-h\nabla\F_3)^{-1}$ exists,
the above
yields $x_{\infty} = J_{h \partial \F_2}\circ(I -h \nabla \F_3)(x_\infty)$.
By the definition
of resolvent in Eq. \eqref{resolvent} this is equivalent to
$\nabla(\F_2 + \F_3)(x_\infty) = 0$. Thus,    the iterates
of this algorithm generates critical points of the underlying ODE.

\section{A classical view on deterministic optimization}\label{sec:classview}

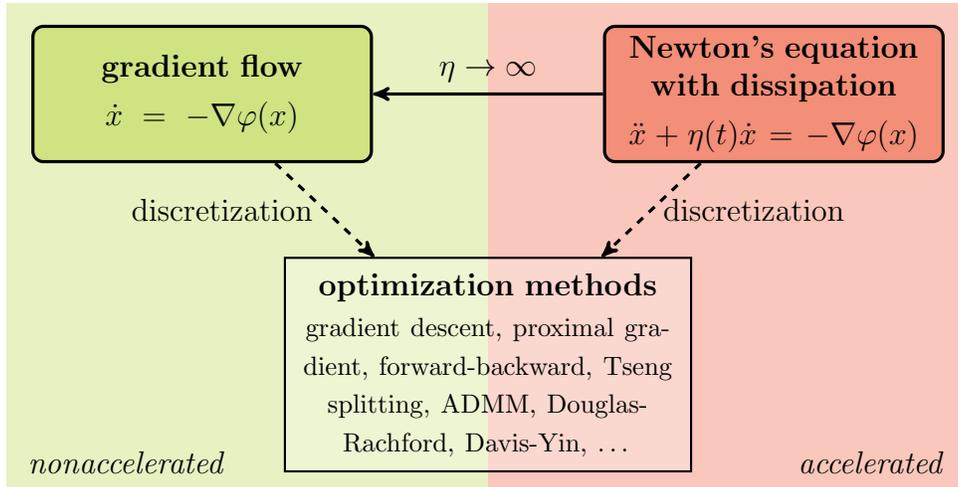
\begin{figure}

\centering

\tikzset{
    >=stealth',
    punkt/.style={
           rectangle,
           rounded corners,
           draw=black, very thick,
           text width=4.5cm,
           minimum height=1.8cm,
           text centered},
    punkt2/.style={
           rectangle,
           %squared corners,
           draw=black, thick,
           text width=5.0cm,
           minimum height=2.8cm,
           text centered
           },
    punkt3/.style={
           rectangle,
           %squared corners,
           %draw=black, thick,
           text width=5.0cm,
           minimum height=2.8cm,
           },
    pil/.style={
           ->,
           thick,
           shorten <=2pt,
           shorten >=2pt,}
}

\definecolor{quad4}{RGB}{30,171,154}
\definecolor{quad1}{RGB}{222,0,109}
\definecolor{quad3}{RGB}{175,209,48}
\definecolor{quad2}{RGB}{235,67,30}

\begin{tikzpicture}[node distance=1cm, auto,]

\fill[quad2!30] (0,-2.9) rectangle (6.4cm,3.6cm);
\fill[quad3!30] (0,-2.9) rectangle (-6.4cm,3.6cm);

\node[punkt3,inner sep=0.2cm,fill=white,opacity=.3] (opt) at (0,-1.2) {};
\node[punkt2,inner sep=0.2cm] (opt) at (0,-1.2)
{{\bf optimization methods}\\[.1em]
\footnotesize{gradient descent, proximal gradient, forward-backward, Tseng splitting, ADMM,
Douglas-Rachford, Davis-Yin, \ldots}};

\node[punkt,inner sep=0pt,fill=quad2!60] (agf) at (3.8,2.4)
    {\bf Newton's equation \\ with dissipation\\[.3em]
         $\ddot{x} + \eta(t)\dot{x} = -\nabla \F(x)$};

\node[punkt,inner sep=0pt,fill=quad3!60] (gf) at (-3.8,2.4)
    {\bf gradient flow \\[.3em] $\dot{x} = -\nabla \F(x)$};

\draw[->,very thick,label] (agf) -- node[above] {$\eta \to \infty$} (gf);
\draw[->,dashed,very thick] (agf) -- node[right] {discretization} (opt);
\draw[->,dashed,very thick] (gf) -- node[left] {discretization} (opt);

\node[] (nonacc) at (-4.8,-2.5) {\emph{nonaccelerated}};
\node[] (acc) at (+5.1,-2.5) {\emph{accelerated}};

\end{tikzpicture}
\caption{Several well-known optimization methods arise
as dicretizations of the gradient flow. Moreover, ``accelerated extensions''
of these methods arise as discretizations of Newton's equation with
a dissipative term; the gradient flow corresponds to a large friction
limit of  the latter.  Therefore, these optimization methods consist
of actual \emph{simulations} of a classical dissipative physical system.
}
\label{overview1}
\end{figure}

At this stage,
let us summarize and interpret the previous results from a physics
perspective. They can be
summarized by Fig.~\ref{overview1}.
Several well-known optimization
methods---including gradient- or proximal-based---are simply different discretizations
of the gradient flow. This includes gradient descent
(see the Appendix), forward-backward, Tseng's splitting,
ADMM, Douglas-Rachford, Davis-Yin, and potentially many other methods.
However, accelerated extensions of these methods arise
as discretizations of Newton's equation with a dissipative term.
Besides the accelerated  proximal methods we introduced,
we also have
Nesterov and heavy ball
(see the Appendix and Ref. \cite{Franca:2019}),
and potentially many others as well.
All these methods are provably first-order
integrators,
i.e., they have a local error $\bigO(h^2)$,
and preserve critical points of the associated ODE.
In addition, let us mention that
we have recently
generalized symplectic integrators
to general dissipative Hamiltonian
systems \cite{Franca:2020}, which  may offer a
systematic approach to construct gradient-based optimization methods.
Therefore, all these optimization methods can be seen as actual
\emph{simulations} of a classical dissipative system;
when someone solves an optimization problem by implementing
one of these methods in a---presumably classical---computer,
there is a simulation of one classical system, i.e.,
dissipative Newton's equation, by another classical
system, i.e., the computer.

There are many interesting open questions even at this
classical level. For instance, given a potential
$\F : \mathbb{R}^n \to \mathbb{R}$ where certain
properties of its landscape are known,
is there a
lower bound on how fast a dissipative system can
approach the ground state?
If so, which specific system would achieve such rate of
convergence?
Given $\F$, what is the best way to dissipate energy, i.e.,
the best damping $\eta(t)$?
What is the tradeoff between energy dissipation and stability?
Answering these type of questions would allow us to
design ``optimal'' optimization algorithms through suitable discretizations
of dissipative physical systems.

A dissipative Newtonian system is deterministic, and so are
the optimization methods obtained as its discretizations.
Such a deterministic approach is suited to convex problems which
have a unique global minimum, or to obtain local minima in the neighborhood of
the initial state.  However, to escape poor local minima in more
complex landscapes, i.e., to solve
more challenging nonconvex problems,
we need to introduce some sort of perturbation or  noise.
Thus, in the next section,
we will generalize these methods to stochastic optimization
settings.  This will allow us to make connections
with Langevin and Fokker-Planck equations
which are ubiquitous in nonequilibrium statistical mechanics.
It is worth noticing that such an approach is in some sense more closely
related  to sampling  than to pure or deterministic optimization.

\section{Stochastic optimization} \label{sec:stochastic}

\subsection{Stochastic gradient}
One of the motivations behind stochastic optimization
is to lighten the
computational burden in computing full gradients over entire
datasets, which is a bottleneck for high-dimensional problems
with large data.  The basic idea dates back to Robbins and Monro
\cite{Robbins:1951} and nowadays is widely used in machine learning,
especially in training neural networks.
Consider replacing the deterministic
problem \eqref{optimization} by its stochastic counterpart
\begin{equation} \label{minE}
\min_{x \in \mathbb{R}^n} \E_{\omega} [ \ell(x ; \omega) ] ,
\end{equation}
where $\omega$ is a random variable from a sample space $\Omega$.
Specifically, suppose we have training data
$\{\omega_1, \dotsc, \omega_N\}$ so
that $\ell_i(x) \equiv \ell(x ; \omega_i)$ is a random variable.
Numerically,
the above expectation is approximated by the empirical mean,
\begin{equation} \label{sloss}
\overline{\ell}(x) \equiv \dfrac{1}{N} \sum_{i=1}^N \ell_i(x),
\end{equation}
which is exact when $N\to\infty$.  Thus,
instead of computing
$\overline{\nabla \ell}(x) =  \tfrac{1}{N} \sum_{i=1}^N \nabla \ell_i(x)$
that may not be feasible, at each iteration of the algorithm we sample
a ``minibatch'' $\B$, of size $S$, drawn uniformly at random---without
replacement---from an index set $\{ 1, \dotsc, N\}$ and compute
the so-called \emph{stochastic gradient}
\begin{equation} \label{sgrad}
\widetilde{\nabla \ell}(x) \equiv
\dfrac{1}{S} \sum_{i \in \B} \nabla \ell_i(x).
\end{equation}
Note that when $S = N$ this becomes the true gradient
of the empirical loss \eqref{sloss}.
Importantly, when the dataset is very large, i.e.,
$S \ll N$ and $N\to\infty$, the \emph{central limit theorem}
comes into play and
\begin{equation}\label{sgrad_noise}
\widetilde{\nabla \ell}(x) = \overline{\nabla \ell}(x) + \xi(x)
\end{equation}
where  $\xi(x) \sim \mathcal{N}\big(0, \Sigma(x)\big)$.
Thus, the stochastic gradient is an
unbiased estimator of the true gradient of the empirical loss.
It is reasonable to assume that the covariance matrix takes the form
\begin{equation} \label{Sigma_x}
\Sigma(x) = \dfrac{1}{S} C(x) C^T(x)
\end{equation}
for some
matrix $C(x)$. We do not know the specific form of $C(x)$, which is data
dependent, however in  principle it can be estimated.
When $N,S \to \infty$, but with the ratio
$S/N \ll 1$ kept fixed, we have $\Sigma \to 0$ and the stochastic gradient
becomes the gradient of the expectation in problem \eqref{minE}---this can
be seen as a thermodynamic limit.

The stochastic gradient \eqref{sgrad} can be implemented into the previous
algorithms quite  easily.
Since $\F_3$ is the only function assumed to be differentiable
in problem \eqref{optimization2}, we consider
\begin{equation} \label{optimization3}
\min_{x \in \mathbb{R}^n} \,
\F_1(x) + \F_2(x) + \E_\omega [ \F_3(x; \omega) ] .
\end{equation}
The entire family of Algorithms~\ref{agenadmmrebal}, \ref{ady} and
\ref{atseng} can be adapted to such case by adding two simple
steps at each iteration, i.e., in the very first line of
the ``for'' loop:
\begin{enumerate}
\item Sample a minibatch $\B \subset \{1, \dotsc, N\}$ of  size $S$ uniformly
at random and without replacement;
\item Replace $\nabla \F_3 \to \widetilde{\nabla \F}_3$ in the
subsequent updates.
\end{enumerate}

\subsection{Stochastic proximal operator}

We can use similar ideas for proximal operators.
As before, in each iteration of the algorithm  we sample a minibatch
$\B$ and define
\begin{equation}
\widetilde{\ell}(x) \equiv \dfrac{1}{S}\sum_{i\in\B} \ell_i(x), \qquad
\ell_i \equiv \ell(x, \omega_i).
\end{equation}
Thus,  at  each iteration, the algorithm has
access to a random function $\widetilde{\ell}(x)$ that presumably
``mimics''  $\E_\omega [ \ell(x; \omega) ]$.
We replace the proximal operator of the empirical loss,
$\prox_{h \overline{\ell}}(x)$, by
its stochastic counterpart
\begin{equation} \label{sprox}
\begin{split}
\widetilde{\prox}_{h \ell}(x) &\equiv \prox_{h \widetilde{\ell}}(x)  \\
&= \argmin_{y}\left( \widetilde{\ell}(y) - \dfrac{1}{2 h}\| y - x\|^2\right) .
\end{split}
\end{equation}
Suppose we introduce stochasticity through $\F_2$ in problem
\eqref{optimization2}, i.e.,
\begin{equation} \label{optimization4}
\min_{x \in \mathbb{R}^n} \,
\F_1(x) + \E_\omega [ \F_2(x; \omega) ] + \F_3(x) .
\end{equation}
Then the family of Algorithms~\ref{agenadmmrebal}, \ref{ady} and
\ref{atseng} are adapted
by adding the following instructions at each iteration:
\begin{enumerate}
\item Sample a minibatch $\B \subset \{1, \dotsc, S\}$ of size $S$ uniformly
at random and without replacement;
\item Replace $\prox_{h\F_2} \to \widetilde{\prox}_{h\F_2}$ in the
next updates.
\end{enumerate}

Note that, also in this case, a similar relation to \eqref{sgrad_noise} holds.
Indeed, from \eqref{proxappro} we get
\begin{equation} \label{sprox_noise}
\begin{split}
\widetilde{ \prox }_{h \ell}(x) &\approx x - h \widetilde{\nabla \ell}(x)\\
&= x - h \overline{\nabla \ell}(x) + h \xi(x) \\
&\approx \prox_{h \overline{\ell}}(x) + h \xi(x).
\end{split}
\end{equation}

\subsection{Langevin and Fokker-Planck equations}

In the deterministic case we have the situation depicted in
Fig.~\ref{overview1}.
In light of the discussion above, introducing a stochastic gradient or
a stochastic proximal operator into these methods is equivalent
to introducing a random perturbation in the associated ODEs.
Thus, the only difference compared to the deterministic case
is that $\nabla \F(x)$ is replaced by a ``stochastic
gradient,'' $\widetilde{\nabla \F}(x)$, during a time interval of one
step size $h$.
As a consequence of Eqs. \eqref{sgrad_noise}
and \eqref{sprox_noise}, we can describe this process by a Brownian motion
provided we account for the correct power of the step size $h$ when discretizing
the noise term. Therefore, we must choose
\begin{equation} \label{disc_noise}
\sqrt{\dfrac{h}{S}} C(x) dW \to \sqrt{\dfrac{h}{S}} C(x_k) \sqrt{h} \,
\epsilon_k \equiv h \xi(x_k)
\end{equation}
where $W$ is a standard Wiener process,
$\epsilon_k \sim \mathcal{N}(0,I)$,
and $\xi(x)$ is the noise
of the stochastic gradient \eqref{sgrad_noise}, or the noise of the stochastic
proximal operator \eqref{sprox_noise}.
The gradient flow \eqref{gradflow} is  thus
replaced by the \emph{overdamped Langevin} equation
\begin{equation} \label{overLang}
dx = - \nabla \F(x) dt + \sqrt{\dfrac{h}{S}} C(x) dW.
\end{equation}
Similarly, the accelerated gradient flow \eqref{secode} is replaced by
the \emph{underdamped Langevin} equation
\begin{subequations} \label{underLang}
\begin{align}
dx &= p dt, \\
dp &= - \nabla \F(x) dt - \eta(t) p dt  + \sqrt{\dfrac{h}{S}} C(x) dW.
\end{align}
\end{subequations}

There is one subtle point about these SDEs.
They have a \emph{multiplicative} white noise which is often ambiguous,
e.g., the It\^o-Stratonovich dilemma.
In our context, it should be noted that the stochastic versions of the
proximal methods previously discussed can be obtained from these SDEs as long
as we discretize the noise consistently with the
gradient of either $\nabla \F_3$ or $\nabla \F_2$, i.e.,
the previous splitting schemes must
be followed carefully by discretizing $C(x)$ appropriately so that
we can combine the noise into the
stochastic gradient $\widetilde{\nabla\F_3}$---in the
case of \eqref{optimization3}---or
the stochastic proximal operator
$\widetilde{\prox}_{h \F_2}$---in the case of \eqref{optimization4}.
Naturally, there is no ambiguity if we assume \emph{additive}
noise, i.e., constant $C$,
which should already provide insights into
these methods, at least qualitatively.

We should also point out that stochastic versions of gradient descent,
Polyak's heavy ball,
and Nesterov---which are gradient-based methods and widely used
in machine learning---also follow
from this  approach.
%, i.e. they
%are discretizations of either \eqref{overLang} or \eqref{underLang}.
In these cases the discretization
is explicit, i.e.,
one has a single gradient, no splitting, and no proximal operators.
Moreover, the
noise term is discretized in the It\^ o sense.
For instance,
stochastic gradient descent is simply an Euler-Maruyama
discretization of \eqref{overLang}.   The stochastic version
of Nesterov arises similarly from \eqref{underLang} by using
\eqref{eulersec}. For heavy ball, one should
use a ``conformal symplectic integrator'' for the
deterministic part of \eqref{underLang}---see Ref. \cite{Franca:2019}---and
compose with an It\^ o discretization of the noise; %
%We mention these details because gradient-based methods are widely
%used in machine learning, especially in training neural networks,
%and they are also related to either \eqref{overLang}
%or \eqref{underLang}---
we provide these
derivations in the Appendix for completeness.

\begin{figure}[t]
\centering
\tikzset{
    >=stealth',
    punkt/.style={
           rectangle,
           rounded corners,
           draw=black, very thick,
           text width=4.2cm,
           minimum height=1.3cm,
           text centered,
           %text=white
           },
    punkt2/.style={
           rectangle,
           %squared corners,
           draw=black, thick,
           text width=5.8cm,
           minimum height=3.2cm,
           text centered
           },
    punkt3/.style={
           rectangle,
           %squared corners,
           %draw=black, thick,
           text width=5.8cm,
           minimum height=3.2cm,
           },
    pil/.style={
           ->,
           thick,
           shorten <=2pt,
           shorten >=2pt,}
}

\definecolor{quad4}{RGB}{30,171,154}
\definecolor{quad1}{RGB}{222,0,109}
\definecolor{quad3}{RGB}{175,209,48}
\definecolor{quad2}{RGB}{235,67,30}

\begin{tikzpicture}[node distance=1cm, auto,]

\fill[quad2!30] (0,0) rectangle (7.5cm,4cm);
\fill[quad3!30] (0,0) rectangle (7.5cm,-4cm);
\fill[quad1!30] (0,0) rectangle (-7.5cm,4cm);
\fill[quad4!30] (0,0) rectangle (-7.5cm,-4cm);

\node[punkt3,inner sep=0.2cm,fill=white,opacity=.3] (opt) at (0,0) {};
\node[punkt2,inner sep=0.2cm] (opt) at (0,0) {
{\bf optimization methods}\\[.1em]
{\footnotesize
(``stochastic'' and/or ``accelerated'') \\
gradient descent,
proximal gradient, \\
forward-backward,
Tseng splitting,
Douglas-Rachford,
ADMM,
Davis-Yin, heavy ball, Nesterov, \ldots
}
};

\node[punkt,inner sep=0pt,fill=quad1!60] (ul) at (-5,3)
    {\bf underdamped Langevin };
\node[punkt,inner sep=0pt,fill=quad4!60] (ol) at (-5,-3)
    {\bf overdamped Langevin };
\node[punkt,inner sep=0pt,fill=quad2!60] (agf) at (5,3)
    {\bf Newton's equation \\ with dissipation };
\node[punkt,inner sep=0pt,fill=quad3!60] (gf) at (5,-3)
    {\bf gradient flow };

\node[inner sep=5pt] (eta) at ($(ul)!0.5!(ol)$) {$\eta \to \infty$};
\draw[-,very thick] (ul) -- (eta);
\draw[->,very thick] (eta) -- (ol);
\node[inner sep=5pt] (eta2) at ($(agf)!0.5!(gf)$) {$\eta \to \infty$};
\draw[-,very thick] (agf) -- (eta2);
\draw[->,very thick] (eta2) -- (gf);
\node[inner sep=5pt] (T) at ($(ul)!0.5!(agf)$) {$ S \to \infty$};
\draw[-,very thick] (ul) -- (T);
\draw[->,very thick] (T) -- (agf);
\node[inner sep=5pt] (T2) at ($(ol)!0.5!(gf)$) {$ S \to \infty$};
\draw[-,very thick] (ol) -- (T2);
\draw[->,very thick] (T2) -- (gf);

\draw[->,dashed,very thick] (ul) -- (opt.north west);
\draw[->,dashed,very thick] (agf) -- (opt.north east);
\draw[->,dashed,very thick] (ol) -- (opt.south west);
\draw[->,dashed,very thick] (gf) -- (opt.south east);

\coordinate (a) at ($(ul) + (-2.5cm,1.3cm)$);
\coordinate (b) at ($(ul)!0.5!(agf)+(-0.00cm,1.3cm)$);
\coordinate (c) at ($(ul)!0.5!(agf)+(+0.00cm,1.3cm)$);
\coordinate (d) at ($(agf) + (+2.5cm,1.3cm)$);
\draw[|-|] (a)--(b) node[midway,above] {\emph{stochastic}};
\draw[|-|] (c)--(d) node[midway,above] {\emph{deterministic}};
\coordinate (a) at ($(ul) + (-2.8cm,1.0cm)$);
\coordinate (b) at ($(ul)!0.5!(ol)+(-2.8cm,+0.00cm)$);
\coordinate (c) at ($(ul)!0.5!(ol)+(-2.8cm,-0.00cm)$);
\coordinate (d) at ($(ol)+(-2.8cm,-1.0cm)$);
\draw[|-|] (a)--(b) node[midway,above,rotate=90] {\emph{accelerated}};
\draw[|-|] (c)--(d) node[midway,above,rotate=90] {\emph{nonaccelerated}};

\end{tikzpicture}
\caption{\label{overview}
Different optimization algorithms arise from different discretizations
(dashed lines) of the \emph{same physical system}.
The left column of the diagram represents
\emph{stochastic processes} described by a Langevin equation,
while the right column represents \emph{deterministic processes} from
(dissipative) classical mechanics---the transition between these phases is controlled
by the ``temperature'' $T \sim 1/S$ where $S$ is the batch size.
The upper row of the diagram corresponds to an underdamped or \emph{accelerated
regime}, while the lower row corresponds to an \emph{overdamped regime} where
acceleration is negligible---the
transition is controlled by the damping coefficient $\eta$.
Discretizations from each colored quadrant yield (variants of) optimization
algorithms in these ``different phases.''
}
\end{figure}
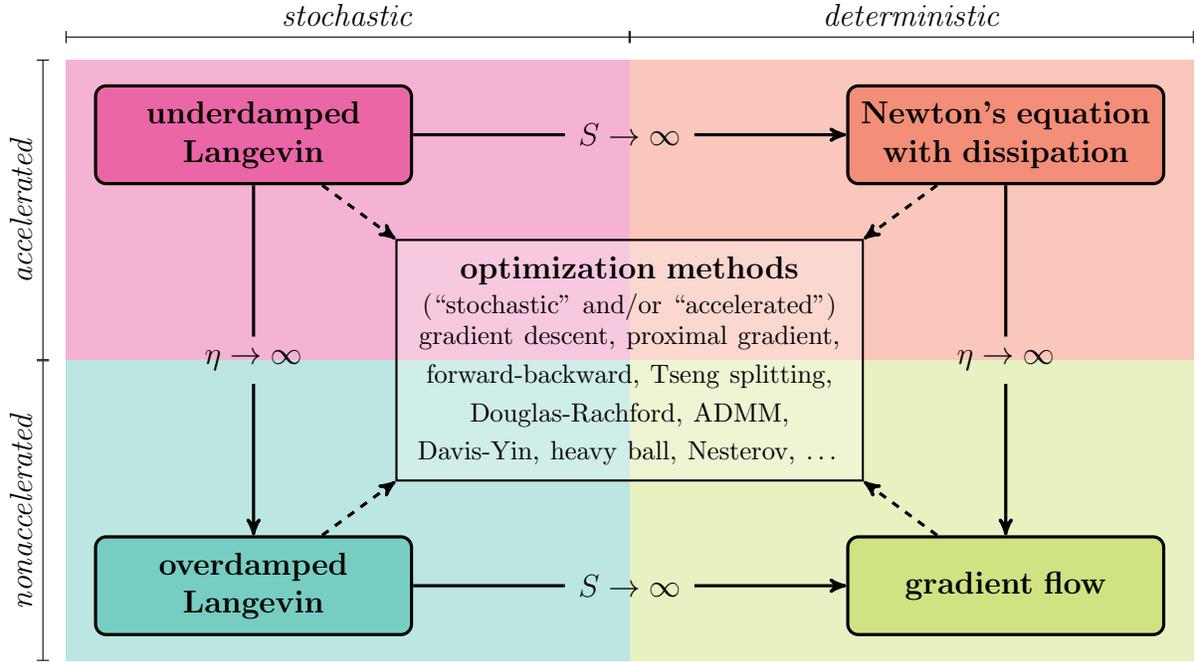

An interesting aspect of the above SDEs
is that the ratio $h / S$ plays the role of
an ``effective temperature'' $T$. This is intuitive: small $S$ means
more noise in the stochastic gradient approximation, which is equivalent
to raising the temperature of the heat bath, while
increasing the step size simply amplifies the noise.
The limit $S \to \infty$---which
presumes $N\to\infty$ and $S/N \ll 1$---corresponds to
removing the heat bath so that
the previous SDEs together with their discretizations
become deterministic.
In similar vein as discussed in Sec.~\ref{sec:high_friction},
the Langevin equation \eqref{overLang} can be recovered
from  \eqref{underLang} in the large friction limit.
Therefore, the overall picture relating all possible variants of these
optimization  methods is depicted in Fig.~\ref{overview};  depending from
which ``quadrant'' one chooses to discretize, and depending which
discretization scheme is chosen, one can obtain an optimization
algorithm with qualifiers
such as ``accelerated,'' ``stochastic,'' or both, or none.
The underdamped Langevin (top left quadrant) is the most
general model that unifies all methods.

Now, we can readily write down the \emph{Fokker-Planck} equation
associated to the above SDEs, which describe
the probability density $P(x, t)$---or $P(x,p,t)$ for the
accelerated methods---of the stochastic process.
From the Chapman-Kolmogorov equation, through a standard derivation, we find
that in the case of \eqref{overLang} we have
\begin{equation} \label{FK1}
\dfrac{\partial P}{\partial t} = \nabla \cdot ( P \nabla \F(x)  )
+ \Delta_{D} P ,
\end{equation}
where  we defined the diffusion matrix and the ``stochastic Laplacian'' as
\begin{equation} \label{Lapla}
D \equiv \dfrac{h}{2S} C(x) C(x)^T, \qquad
\Delta_D P \equiv \sum_{i,j} \dfrac{\partial^2}{\partial x_i \partial x_j}
(D_{ij} P).
\end{equation}
Similarly, the Fokker-Planck equation associated to the underdamped
Langevin \eqref{underLang} is given by
\begin{equation} \label{FK2}
\dfrac{\partial P}{\partial t} = - \nabla_x \cdot \left( p P \right)
+ \nabla_p\cdot\left( P \nabla \F(x)   + \eta(t) p P \right)
+ \Delta_{D} P .
\end{equation}
However, here
$\Delta_D P = \sum_{ij} D_{ij}(x) \partial_{p_i}\partial_{p_j} P$ since
the noise is coupled only to the momenta.

In Eqs. \eqref{overLang} and \eqref{underLang} the Brownian motion
arises from a stochastic gradient or stochastic proximal operator
approximation,
hence the diffusion dependence on the minibatch and step sizes.
Alternatively, one can directly perturb the
deterministic optimization methods %---i.e., without  any  gradient
%approximation---
with a Gaussian noise. Both situations are conceptually similar except that
in the latter case
we have discretizations of the usual Langevin dynamics
where the Brownian motion is controlled independently  by a heat bath.
The noise term of Eq. \eqref{underLang}
is thus  replaced by the standard
$\sqrt{2 D} dW$, which obeys Einstein's fluctuation-dissipation
relation $D = \eta \, k_B T$, assuming $\eta=\mbox{const.}$
The overdamped limit yields Eq. \eqref{overLang} with the same term
but $D = k_B T$. In this setting, we have ``sampling
methods'' where the dynamics first relax to equilibrium
and then perform random excursions around a local minimum of $\F$.
From this perspective, the stochastic
optimization methods are conceptually closer to a sampling approach, however
with smaller and limited control over the noise.

\section{Numerical Experiments}
\label{sec:numerical}

\subsection{Order of accuracy}
\begin{figure}  \centering
\includegraphics[width=0.5\textwidth]{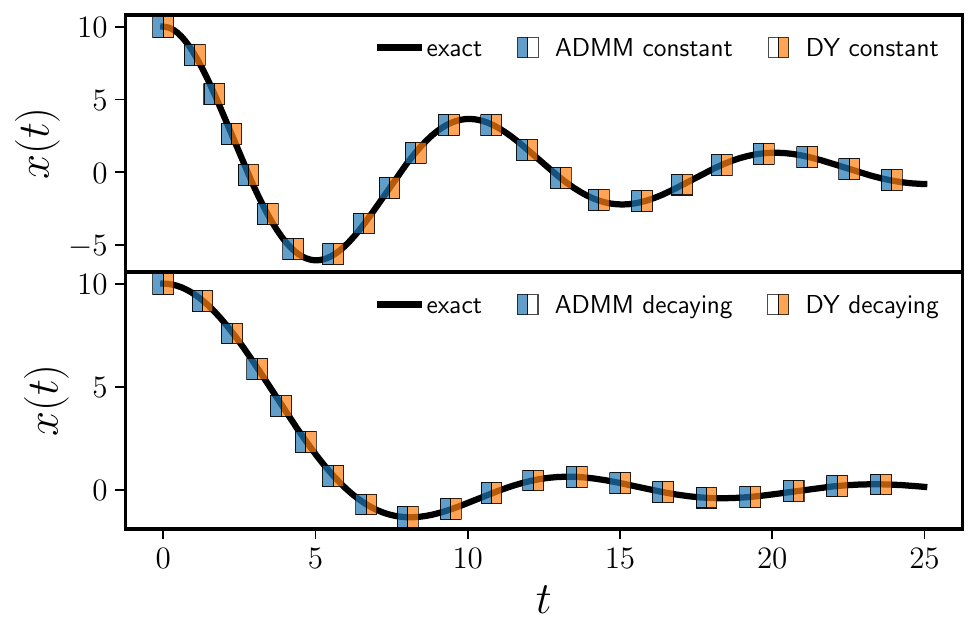}%
\includegraphics[width=.5\textwidth]{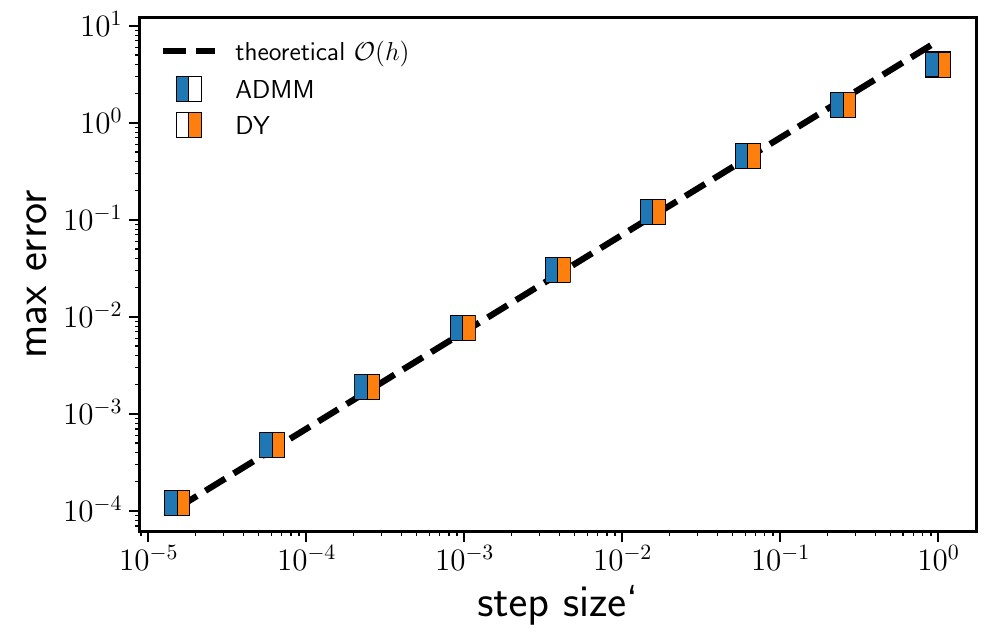}
\put(-440,100){\emph{(a)}}
\put(-440,35){\emph{(b)}}
\put(-180,80){\emph{(c)}}
\caption{%
Comparison of the discretization provided by
Algorithms~\ref{agenadmmrebal} and \ref{ady} versus the exact solution of
the second-order gradient flow \eqref{secode} with \eqref{compos_quad}.
\emph{(a)} Constant damping with $\eta = 0.2$. %; exact solution
%is \eqref{sol_agf_const}.
\emph{(b)} Decaying damping with $r=3$.
%exact solution is \eqref{sol_agf_decaying}.
In both cases we chose step size $h=0.01$.
\emph{(c)}
Maximum error \eqref{max_error}
over a range of step sizes. We consider constant damping
under the  setting of Fig.~\ref{exact_discrete}a---%
similar results hold for decaying damping.
\label{exact_discrete}
\label{quad_error}
}
\end{figure}

We start by illustrating how accurate the previous optimization
methods approximate the gradient flow ODEs.
We consider Algorithms~\ref{agenadmmrebal} and \ref{ady}---recall that
the latter reduces to (accelerated versions of) Douglas-Rachford and
forward-backward
as particular cases, so we are indirectly testing them as well.
We compare these methods with
the exact solutions
\eqref{sol_agf_const} and \eqref{sol_agf_decaying}. To this end, we  choose
a composite function:
\begin{equation} \label{compos_quad}
\F = \F_1 + \F_2 + \F_3, \qquad
\F_i(x) = \omega_i^2  x^2 / 2 ,
\end{equation}
with $\omega_1 = 1/2$, $\omega_2 = 1/3$
and $\omega_3 = 1/5$---recall that there is \emph{splitting} of these
individual functions in the algorithms.
As we can see in Figs.~\ref{exact_discrete}a and \ref{exact_discrete}b,
these methods closely  match the exact solutions.
%Moreover, both methods
%yield the same results, in agreement with our  derivations that although
%both are different discretizations of the same system, they are numerical
%integrators of the same order of accuracy.
Naturally, for large step sizes there might be
significant deviations.
However, since these methods were proven to
be first-order  integrators, they are accurate up to
a global error $\bigO(h)$.
To verify this, for a given step size $h$
we compute
\begin{equation} \label{max_error}
\max_{k\in [0,K]} |x(t_k) - x_k|
\end{equation}
for $K = t_{\textnormal{max}}/h$, where $t_{\textnormal{max}}$ is
a fixed simulation time (we choose again $t_{\textnormal{max}}=25$).
We thus compare this maximum error over the entire history of the system
against a range of  step sizes. Results for the
constant damping case are shown in Fig.~\ref{quad_error}c---%
we performed the same simulation with decaying
damping but the plot is nearly indistinguishable from this one.

In these experiments there is no visible difference
between---accelerated---Davis-Yin and ADMM. This is
also explained by our theoretical results since both methods consist of
numerical integrators of the same order of accuracy and to the same
physical system.
%Let us mention that we also verified similar results for the
%base algorithms, i.e., for the nonaccelerated variants, in relation
%to the first-order gradient flow \eqref{gradflow}---we omit these
%results since they are visually the same as the ones above.

\subsection{Langevin approximation}
\begin{figure}
\centering
\includegraphics[width=0.49\textwidth]{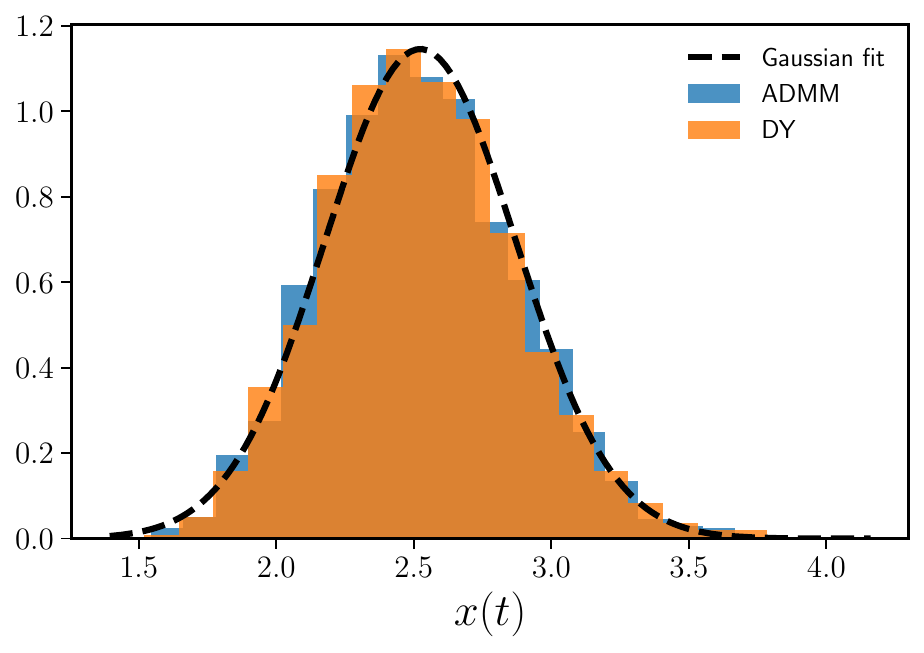}%
\includegraphics[width=0.519\textwidth]{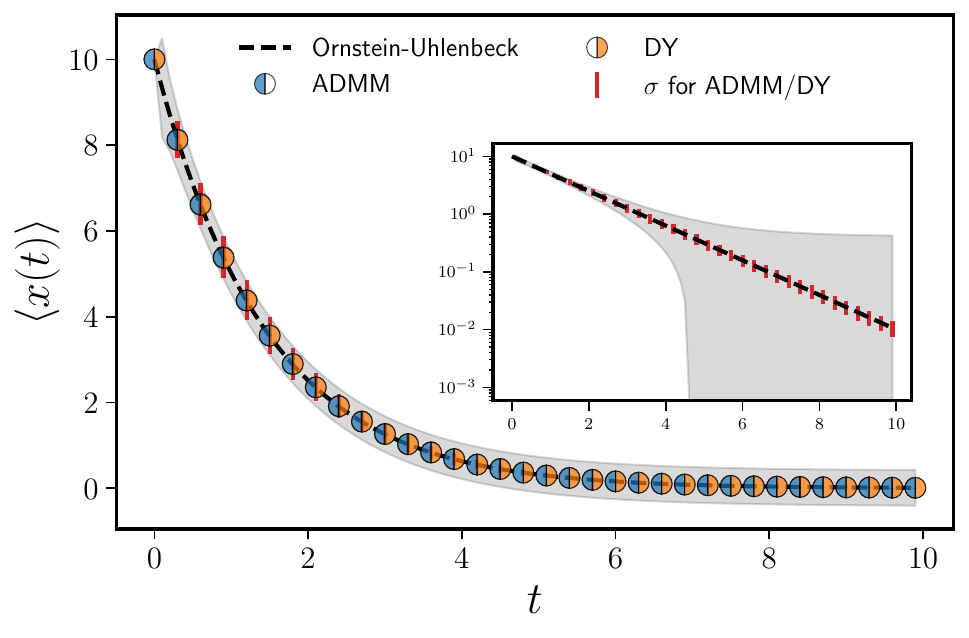}
\put(-450,50){\emph{(a)}}
\put(-200,50){\emph{(b)}}
\caption{\emph{(a)} Histogram of stochastic variants of Algorithms~\ref{agenadmmrebal}
and \ref{ady} (with $\gamma_k = 0$) for a quadratic problem which is
modelled by an Ornstein-Uhlenbeck process; we considered \eqref{F3stoch}
with $N=1000$, minibatch  of size $S=1$,  step size $h=0.1$, initial
position $x_0= 10$ and  $2000$ Monte Carlo runs.
We show the simulation for $t = 2$. The
Gaussian fit has mean $\mu = 2.52$ and variance $\sigma = 0.35$.
The theoretical prediction from Eq. \eqref{ou_sol} is $\mu = 2.51$,
in close agreement.
\emph{(b)}
Same type of simulation  for
$t \in [0,10]$. The dashed line is the mean of the Ornstein-Uhlenbeck
process \eqref{ou_sol}, and the shaded area indicates the standard deviation
with $D=4$.
The markers are the means of ADMM and
Davis-Yin, and the vertical red lines are standard deviations for both
methods, which are very close.  The inset is exactly the same plot but with
markers omitted and the $y$-axis on log scale.
\label{hist}
\label{means}
}
\end{figure}

We now wish to verify the Langevin approximation to the stochastic variants
of the optimization methods introduced in the previous section.
Consider thus a quadratic function as in Eq. \eqref{compos_quad}, again with
$\omega_1 = 1/2$ and $\omega_2  = 1/3$, but now $\F_3$ is generated
randomly so as to obtain a problem in the form \eqref{optimization3}.
We set
\begin{equation}
\label{F3stoch}
\F_3(x) = \dfrac{1}{2N}\sum_{i=1}^N \theta_i^2 x^2
\end{equation}
with $\theta_i$ sampled
uniformly from the interval  $[0,1]$. We choose
$N = 1000$ and $\F_3$ is kept fixed in all
simulations, i.e., we have an ``empirical mean'' as in \eqref{sloss}.
To obtain stochastic gradients \eqref{sgrad},
at each iteration of Algorithms~\ref{agenadmmrebal} and \ref{ady}
we sample a minibatch of $\theta_i$'s
of size $S = 1$---this is the only source of stochasticity.
For simplicity, we focus on the nonaccelerated
methods ($\gamma_k = 0$) which are modelled by the overdamped
Langevin \eqref{overLang}.
Since the function is quadratic, if we assume a constant $D$
into \eqref{FK1}---which we know is not the
case but it should already capture the qualitative behavior---we then
have an Ornstein-Uhlenbeck process whose
probability density, $P(x,t | x_0, t_0=0)$,
is given by
\begin{equation} \label{ou_sol}
\sqrt{ \dfrac{\lambda}{2 \pi D (1-e^{-2\lambda t})}}
\exp\left\{ -\dfrac{\lambda \left(x  - x_0 e^{-\lambda t}\right)^2}
{2 D (1-e^{-2\lambda t})}  \right\},
\end{equation}
with $\lambda = \omega_1^2 + \omega_2^2 +
\tfrac{1}{2N}\sum_{i}^{N}\theta_i^2$.
In our simulations we choose a step  size $h = 0.1$, maximum time
$t_{\textnormal{max}}  = 10$, and initial position
$x_0 = 10$.  We consider $2000$ Monte Carlo runs, with $\F_3$ fixed so that
only the minibatch changes in each simulation.
In Fig.~\ref{hist}a we show histograms for
stochastic ADMM and Davis-Yin
%---again, we are indirectly  testing
%forward-backward and Douglas-Rachford as well---
for a single time
instant ($t=2$). We can see that both
methods give essentially the same results and agree with the
out-of-equilibrium Gaussian
prediction \eqref{ou_sol}.
Similar results hold for other time instants, as illustrated
in Fig.~\ref{means}b. Here we set $D = 4$  to obtain
the shaded gray area---this value was estimated with the simulation
for $t=1$, and we picked this time since the variance was
sufficiently large.  We also indicate the standard  deviation of
both methods (vertical red lines) which actually vary with
$x$, however everything is well-described
by Eq. \eqref{ou_sol}, in agreement with the Langevin approximation.
We mention that simulations with
other values of the batch and step sizes
were also considered, verifying similar results but
where the Gaussian approximation is more or less peaked according to the
the scaling $h/S$ from Eq. \eqref{Lapla}.

\subsection{Machine learning experiments}
We now wish to verify whether the accelerated
methods we introduced are able to
achieve faster convergence
compared to the base methods, which are the actual known methods in
the literature.
According to the discussion in Sec.~\ref{sec:rates}, we expect
that this might be the case.
We focus on two types of damping: constant  \eqref{hbdamp},
and decaying \eqref{nagdamp}.
When nothing is specified it means that no acceleration is used, i.e.,
$\gamma_k=0$ in Eq. \eqref{xhat}, which yields the base methods.
We consider ADMM (Algorithm~\ref{agenadmmrebal}),
DY (Algorithm~\ref{ady}) and Tseng (Algorithm~\ref{atseng}).
Forward-backward (FB) is Algorithm~\ref{ady} with $\F_1 = 0$ and
Douglas-Rachford (DR) with $\F_3 = 0$.

\begin{figure}
\centering
\includegraphics[scale=0.6]{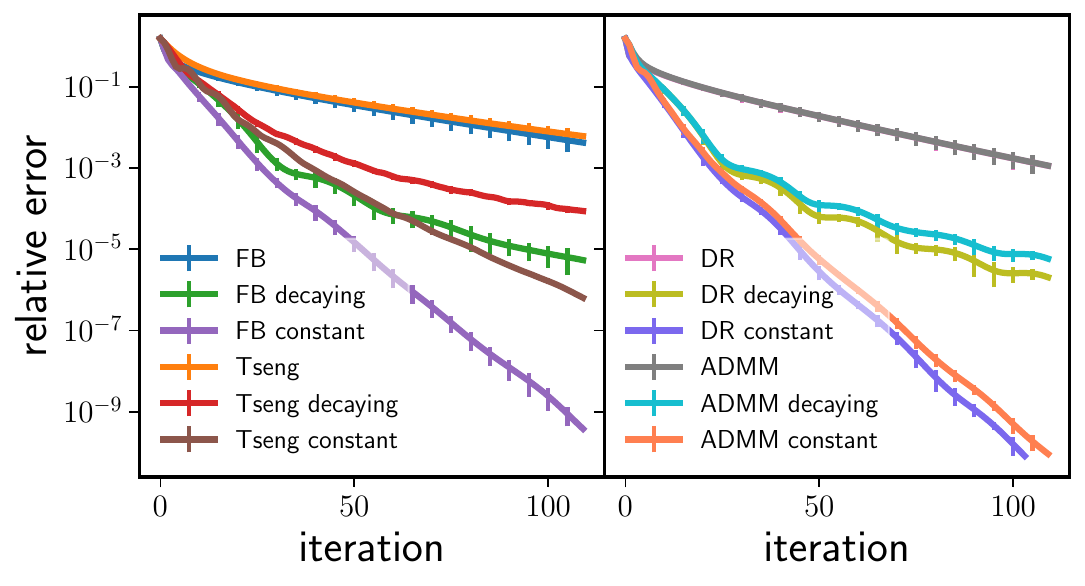}
\caption{\label{lasso}
Performance of several methods on problem~\eqref{lassoprob}.
We perform 10 Monte-Carlo runs and show the mean and standard deviation (error
bars) of
the error
$|\F(x_k) - \F^\star|/\F^\star$ versus the iteration $k$.
%Here $\F$ is the objective function and $\F^\star$ is the minimum obtained
%with an independent and reliable solver (CVXPY).
The accelerated variants have faster convergence.
}
\end{figure}

Consider a \emph{LASSO regression} problem
that is of fundamental importance
to machine learning and statistics:
\begin{equation}\label{lassoprob}
\min_{x\in\mathbb{R}^n} \F(x), \qquad \F(x)
\equiv \dfrac{1}{2}\| A x - b\|^2 + \alpha \| x\|_1 ,
\end{equation}
where $A \in \mathbb{R}^{m\times n}$ is a given matrix, $b \in \mathbb{R}^m$
is a given signal, and $\alpha > 0$ is a coupling constant.
Here $\| \cdot \|_1$
denotes the $\ell_1$-norm, known to induce ``sparsity''---this
function is not differentiable although
its proximal operator has a well-known closed form
solution called \emph{soft thresholding}.
Following Ref. \cite[Sec.~7.1.3]{Parikh:2013},
we generate data by sampling
$A_{ij} \sim \mathcal{N}(0,1)$
and then normalize its columns to have unit norm.
We sample
$x_{\bullet} \in \mathbb{R}^{n}\sim \mathcal{N}(0,1)$ with sparsity
level $95\%$ (only $5\%$ of its entries are nonzero) and then add
noise to obtain the observed signal
$b = A x_{\bullet} + e$ where $e \sim \mathcal{N}(0, 10^{-3})$.
We choose dimensions $m=500$ and $n=2500$---the signal-to-noise ratio
is $\approx 250$, and  $x_\bullet$ has $125$ nonzero entries.
We set
$\alpha= 0.1 \alpha_{\textnormal{max}}$ where
$\alpha_{\textnormal{max}} = \| A^T b \|_{\infty}$ is the largest value
such that \eqref{lassoprob} admits a nontrivial solution---the
factor $0.1$ was verified to yield good results
after few trials.
We evaluate the algorithms by computing the relative error
$|\F(x_k) - \F^\star| / \F^\star$ where
$\F^\star$ is the solution obtained with an independent and
reliable solver---we use the default implementation of CVXPY
which is a standard optimization library in the Python language.
For all algorithms, we choose step size $h = 0.08$, which is the largest
choice such that all algorithms converge.
For decaying damping we choose $r=3$ in \eqref{nagdamp}---this choice
is standard but we verified that other
values did not improve---and for constant damping we choose
$\eta=0.5$  in \eqref{hbdamp}---this value was
tuned with a rough grid search.
In Fig.~\ref{lasso} we report the mean
and standard deviation (small error bars) across 10 randomly
generated instances of
problem~\eqref{lassoprob}.
The figure shows that the accelerated variants of each method improve
over the base method.
In particular, the constant damping variants are the fastest.

Next, we consider a matrix completion problem
which is also of fundamental importance
in machine learning. The goal is to reconstruct a \emph{low-rank} matrix where
we are only allowed to observe a few of its entries.
Moreover, we assume these entries
are constrained to specified range. More precisely, suppose that for
a low-rank matrix $M \in \mathbb{R}^{n\times m}$
we observe entries $(i,j)$ that are collected in a set $\Omega$:
let $\PP:\mathbb{R}^{n\times m} \to \mathbb{R}^{n\times m}$
be the projection onto the support of observed entries.
The observable data matrix is thus $\Mobs = \PP(M)$, where
$\PP(M)_{ij} = M_{ij}$ if $(i,j) \in \Omega$
and $\PP(M)_{ij} = 0$ otherwise.
The goal is to estimate the missing entries of $M$.
This can be done \cite{Cai:2010} by solving the
convex problem
$\min_X  \| X \|_{*} $ such that $\PP(X) = \PP(M)$---here $\| X \|_{*}$
denotes the nuclear norm.
We consider a modification of this approach by imposing
constraints $a \le X_{ij} \le b$ for given constants $a$ and $b$.
Specifically,
\begin{equation}
  \label{matcomp}
  \min_{X\in\mathbb{R}^{n\times m}}  \underbrace{\alpha \| X\|_*}_{\F_1}  +
  \underbrace{\mathcal{I}_{[a,b]}(X)}_{ \F_2 }  +
  \underbrace{\dfrac{1}{2} \|  \PP(X) - \PP(M) \|_{F}^{2}}_{\F_3}
\end{equation}
where $\| \cdot \|_F$ denotes the Frobenius norm,
$\mathcal{I}_{[a,b]}(X) = 0$ if
$ a \le X_{ij} \le b$  and $\infty$ otherwise. A higher  $\alpha > 0$
induces lower rank solutions.
This problem can be solved with Algorithms~\ref{agenadmmrebal} and
\ref{ady} with the proximal operator
\begin{equation}
\prox_{h \| \cdot \|_*}(X) = U D_{h}(\Sigma) V^T ,
\end{equation}
where $X = U \Sigma V^T$ is the singular value decomposition of $X$ and
$D_{h}(\Sigma)_{ii} = \max\{ \Sigma_{ii} - h, 0\}$
(see \cite{Cai:2010} for details).
The proximal operator of $\F_2$ is just the projection
\begin{equation}
\prox_{h \mathcal{I}_{ [a,b] } }(X)_{ij}
= \max\{a, \min(X_{ij}, b)\}.
\end{equation}
Moreover, $\nabla \F_3(X) = \PP(X - M)$. In our methods
we choose the following stopping criterion:
\begin{equation}\label{tol}
  \| X_{k+1} - X_k \|_F   \big/   \| X_k\|_F  \le \epsilon
\end{equation}
where $\epsilon$ is a small tolerance.
To evaluate  performance we report the relative error
\begin{equation}\label{rel_error}
 \| X_k - M \|_F \big/  \| M\|_F .
\end{equation}

\begin{figure}[t!]
\centering
\includegraphics[scale=0.58]{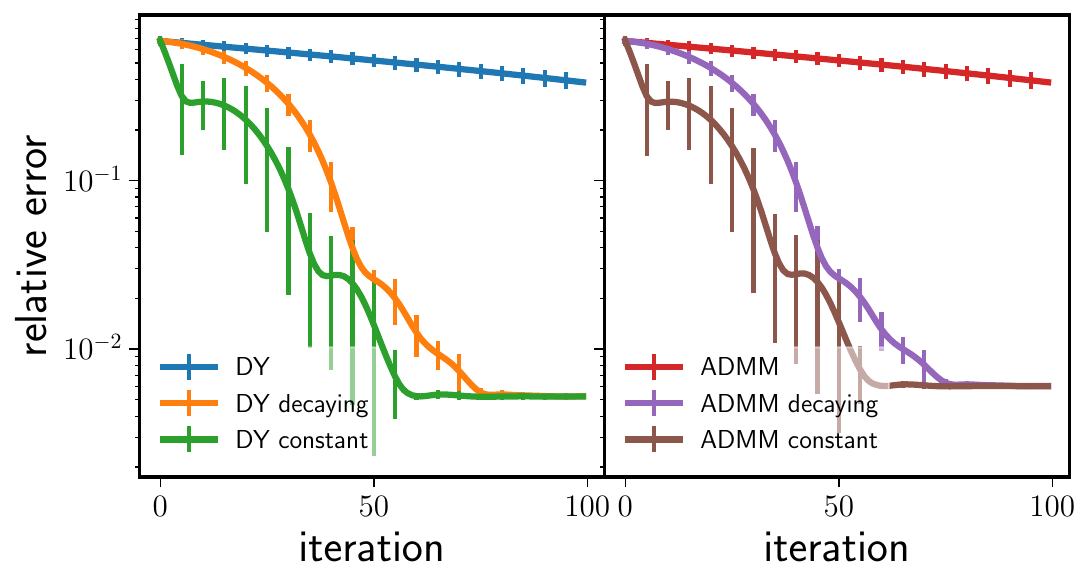}
\caption{\label{mc1}%
Convergence of different algorithms on problem \eqref{matcomp}.
We perform $10$ Monte Carlo runs and
indicate the mean and standard deviation (error bars) for the relative error
\eqref{rel_error}
versus iteration  $k$.
Note the improvement of the accelerated methods.}
%\end{figure}
%
\vspace*{1em}
%\begin{figure}[t]\centering
\includegraphics[width=0.49\textwidth]{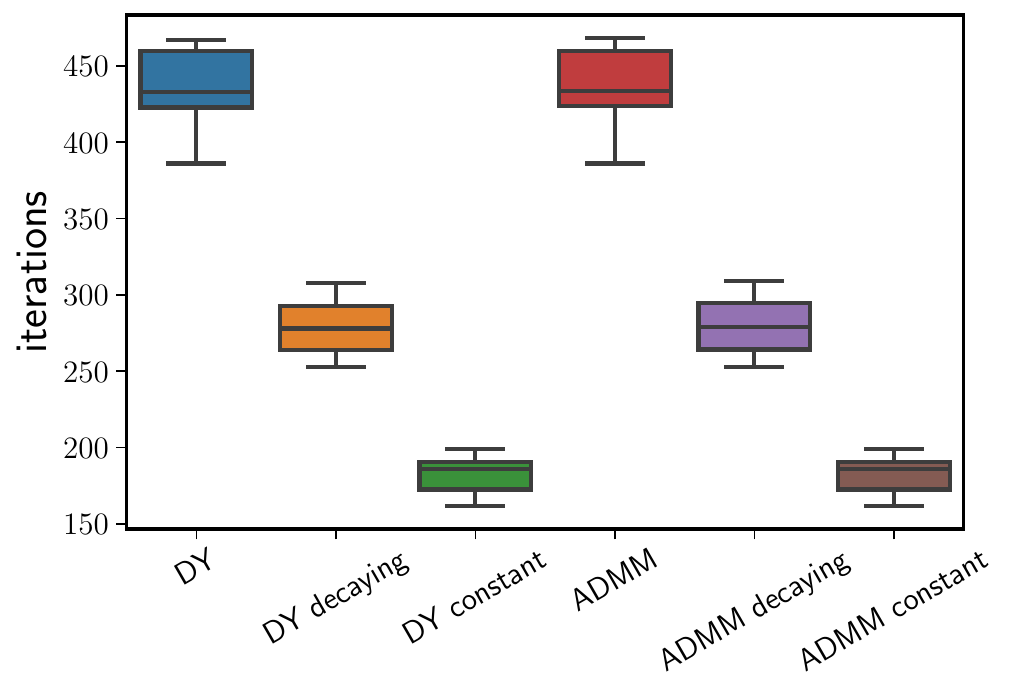}
\includegraphics[width=0.49\textwidth]{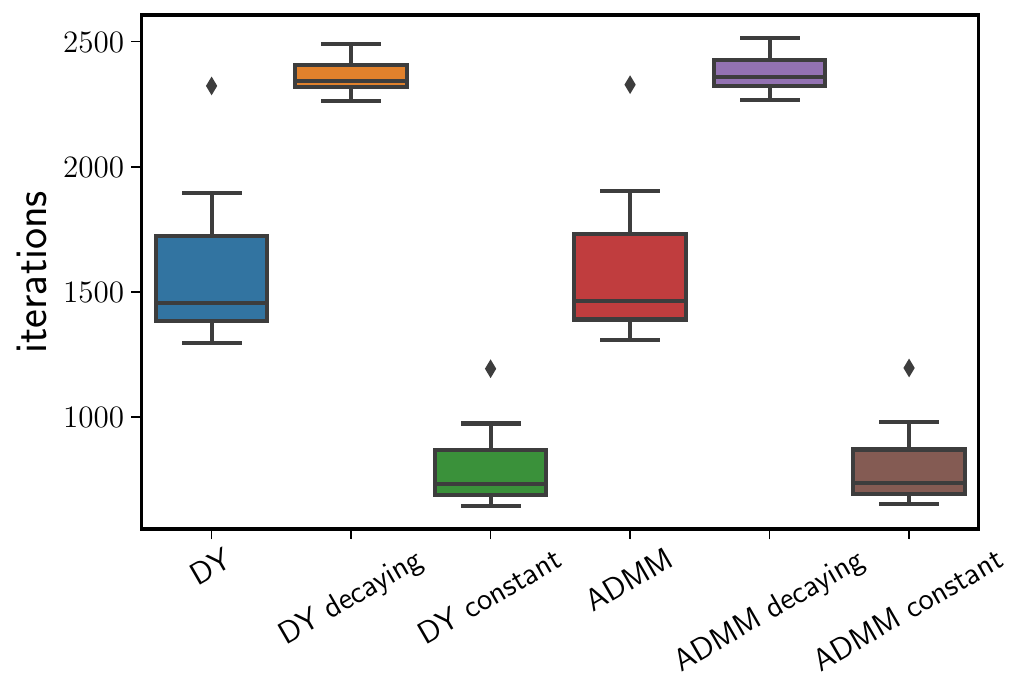}
\put(-470,10){\emph{(a)}}
\put(-220,10){\emph{(b)}}
\caption{
\label{mc12}%
\emph{(a)} Number of iterations needed to reach the
termination tolerance for the problem in Fig.~\ref{mc1}.
\emph{(b)} Number of iterations needed to reach the
termination tolerance for the problem in Fig.~\ref{mc2}.
\label{mc22}
}
%\vspace*{1em}
%\end{figure}
%
%\begin{figure}[t]
%\centering
\includegraphics[scale=0.58]{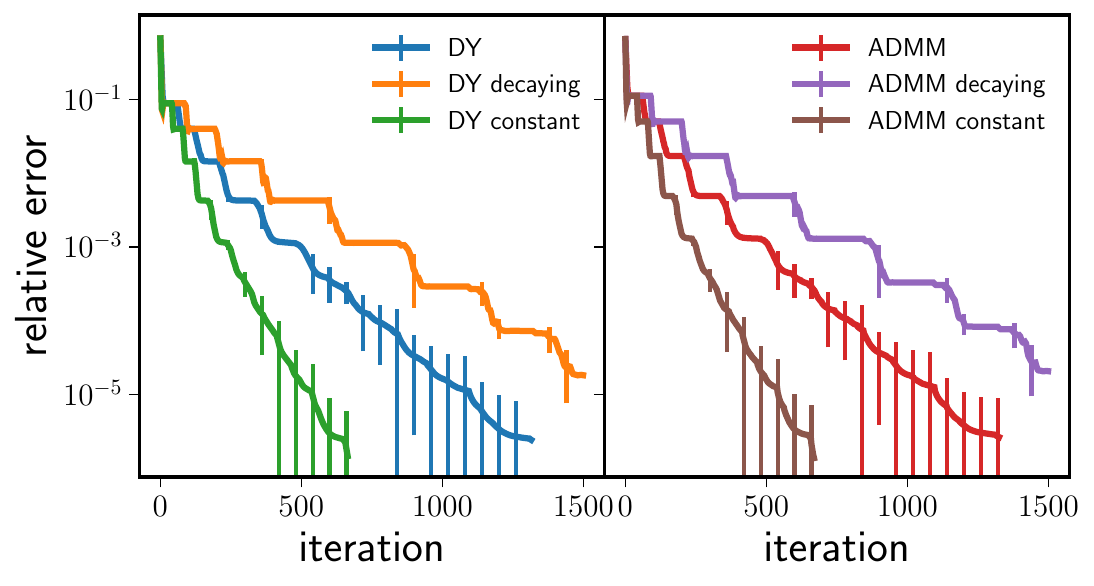}
\caption{Performance of different algorithms on problem~\eqref{matcomp}
%under the same setting as in Fig.~\ref{mc1} but
with  annealing on $\alpha$.
\label{mc2}
}
\end{figure}

Following Ref. \cite{Cai:2010}, we generate a low-rank matrix as $M = L_1 L_2^T$
where $L_1, L_2 \in \mathbb{R}^{100 \times 5}$ with entries sampled
i.i.d. from $\mathcal{N}(3, 1)$.
Thus $M$ has rank $5$ (with probability one) and each entry is positive
with high probability (each test instance was checked to
have positive entries).
We sample $s n^2$ entries of $M$ uniformly, with
a sampling ratio $s=0.4$ (i.e. only $40\%$ of the matrix $M$ is observed).
We choose
\begin{equation}
\begin{split}
a &= \min\{ (M_{\textnormal{obs}})_{ij} \} - \sigma/2 , \\
b &= \max\{ (M_{\textnormal{obs}})_{ij} \} + \sigma/2,
\end{split}
\end{equation}
where $\sigma$ is the
standard deviation of all entries of $M_\textnormal{obs}$.
In terms of algorithm's parameters, we choose
a step size $h=1$ for all methods, $r=3$ for
decaying damping \eqref{nagdamp},
and $\eta=0.1$ for constant damping \eqref{hbdamp}---the choice of step size
is standard for this type of problem with proximal methods, while
the choice of damping was obtained with a rough grid search and provided
good results. In the stopping criterion \eqref{tol} we choose
$\epsilon = 10^{-10}$.
In Fig.~\ref{mc1} we report the mean
and standard deviation (error bars) across 10 randomly generated instances of
problem~\eqref{matcomp} with $\alpha = 3.5$;
all methods terminate successfully and recover a matrix with the correct rank
(for this choice of $\alpha$)
and a final relative error of $\approx 5\times 10^{-3}$.
The number of iterations of each method to achieve the tolerance
error are reported in  Fig.~\ref{mc12}a.

To obtain more accurate solutions
we consider ``annealing'' the parameter $\alpha$.
We wish to
verify if the accelerated methods can still speedup convergence
in this scenario.
We follow the procedure of \cite{Goldfarb:2011} which is as follows.
Given a sequence
$\alpha_1 > \alpha_2 > \dotsm > \alpha_{L} = \bar{\alpha} > 0$
for some $\bar{\alpha}$,
we run each algorithm with $\alpha_j$
and then use its solution as a warm start for the solution to the next run
with $\alpha_{j+1}$ (all other parameters are kept fixed).
Starting with $\alpha_0 = \delta \| \Mobs \|_F$ for some $\delta\in(0,1)$ we
use the schedule
\begin{equation} \label{anneal_alpha}
\alpha_{j+1} = \max\{\delta \alpha_j, \bar{\alpha}\}
\end{equation}
until
reaching $\bar{\alpha}$.
In our tests we choose $\delta = 0.25$ and $\bar{\alpha}=10^{-8}$---these
parameters were tuned with a few trials and proved to yield good results.
The remaining parameters are the same as those used in creating Fig.~\ref{mc1},
except that for the constant damping variants we now use $\eta=0.5$---in this
example, overdamping the system yield better results.
We report the convergence of different methods in Fig.~\ref{mc2}.
%across $10$
%randomly generated instances of this problem.
All methods successfully reach the termination tolerance, as for the
previous test,
but now achieve a better reconstruction accuracy.
The total number of iterations in this case
are shown in Fig.~\ref{mc22}b.
Note that, in this experiment, the decaying damping variants did
not improve over the nonaccelerated method, but the constant damping
variants still provided a significant improvement.
This has to do with the ``hardness'' of the problem, i.e., due to the choice
of rank and sampling ratio the objective function should have
high curvature and this is why the overdamped case is
faster (recall Fig.~\ref{solutions_fig}).

\begin{figure}[t]
\centering
\includegraphics[scale=0.6]{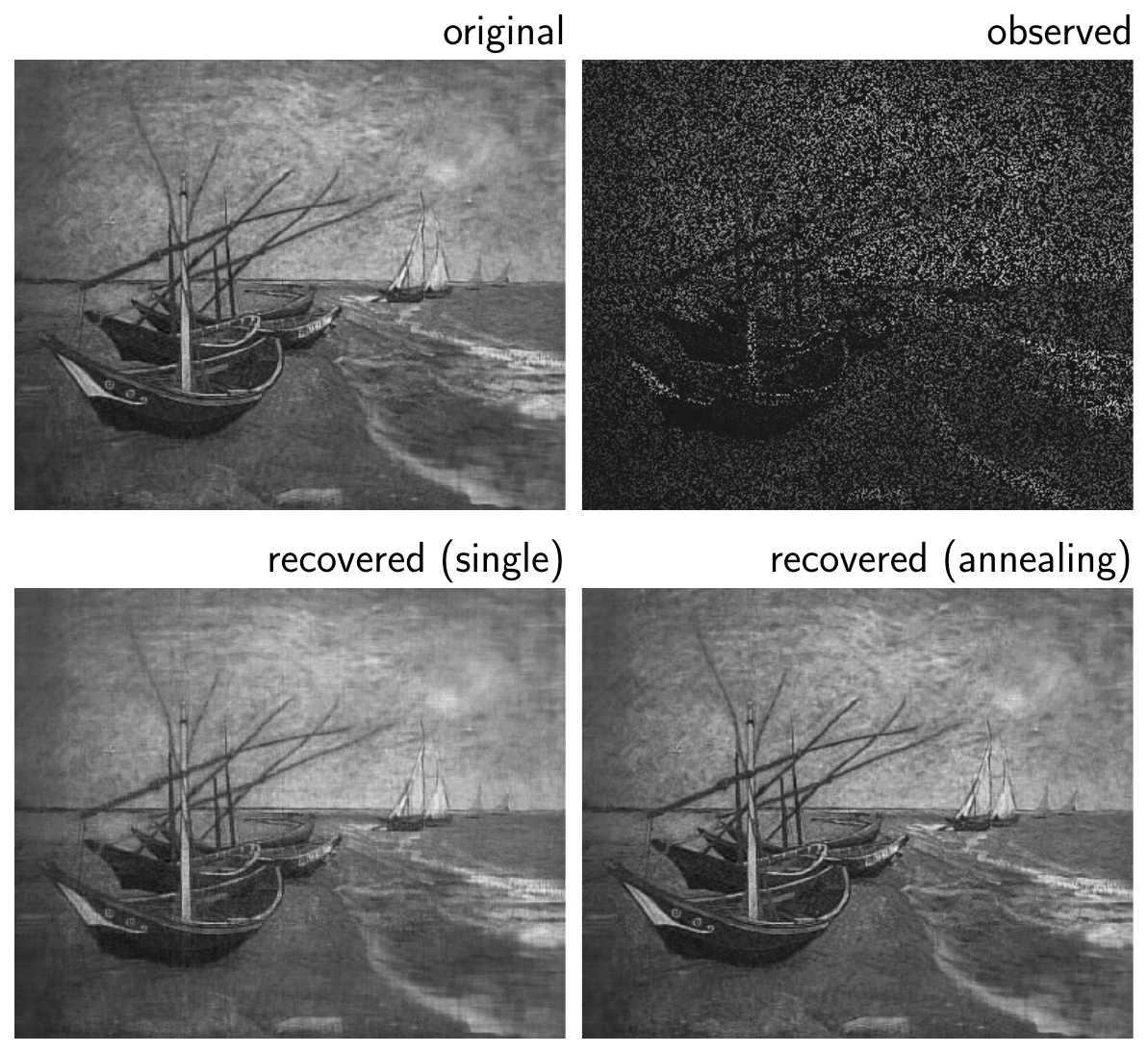}
\put(-330,302){\emph{(a)}}
\put(-160,302){\emph{(b)}}
\put(-330,145){\emph{(c)}}
\put(-160,145){\emph{(d)}}
\caption{%%%
\emph{(a)} Image of $974 \times 1194$ pixels with rank $r=70$.
\emph{(b)} Observed data with sampling ratio $s=0.3$.
\emph{(c)} Reconstructed image by running the previous
algorithms a single time.
The relative error is $8 \times 10^{-2}$---%
see Fig.~\ref{boat_convergence1} for convergence rates.
\emph{(d)} Reconstructed image with annealing.
The relative error is $1.6 \times 10^{-4}$---see
Fig.~\ref{boat_convergence2} for convergence rates.
The correct rank is also recovered by all methods.
\label{boats_fig}
}
\end{figure}

We now provide a real world example where
the goal is to reconstruct a partially observed image.
We pick a gray scale image, and
by a singular value decomposition we truncate its first 70 eigenvalues.
The end result is the image shown in Fig.~\ref{boats_fig}a,
represented by a $974 \times 1194$ matrix with entries between
$0$ and $1$ and of low rank $r=70$---%
this image will not be used as input data
but will serve as the ground truth $M$.  Then, we sample some of its entries
uniformly with sampling ratio $s=0.3$, obtaining
the image shown in Fig.~\ref{boats_fig}b---this is the input
data $M_{\textnormal{obs}}$
and has only $30\%$ of the original entries.
This problem is actually hard, e.g., the number of effective
degrees of freedom is $r(n+m-r)/(s n m ) \approx 0.42$
(see \cite{Cai:2010} for details).

\begin{figure}[b!]\centering
\includegraphics[width=0.5\textwidth]{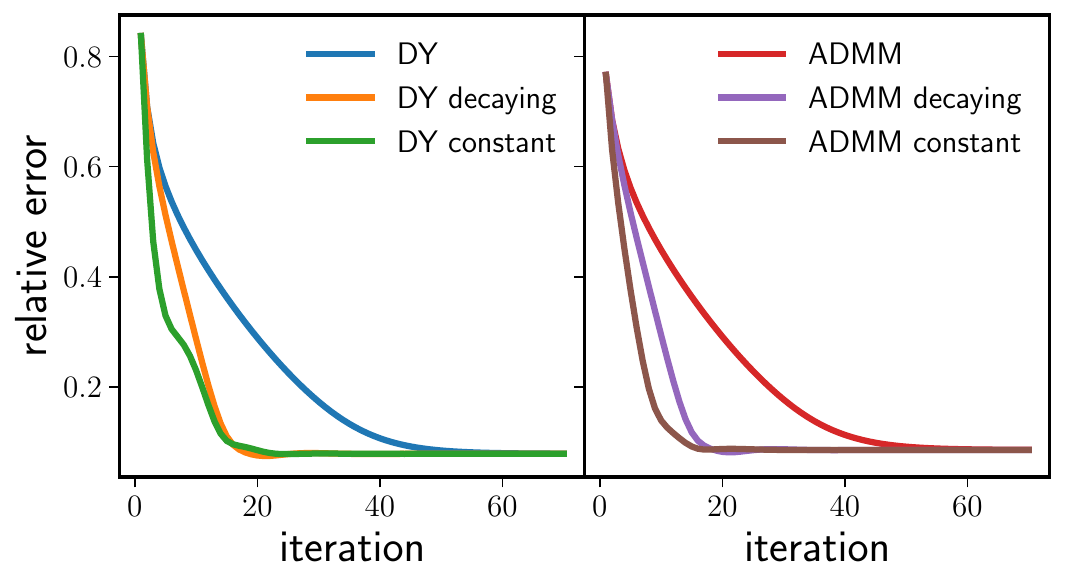}%
\includegraphics[width=0.5\textwidth]{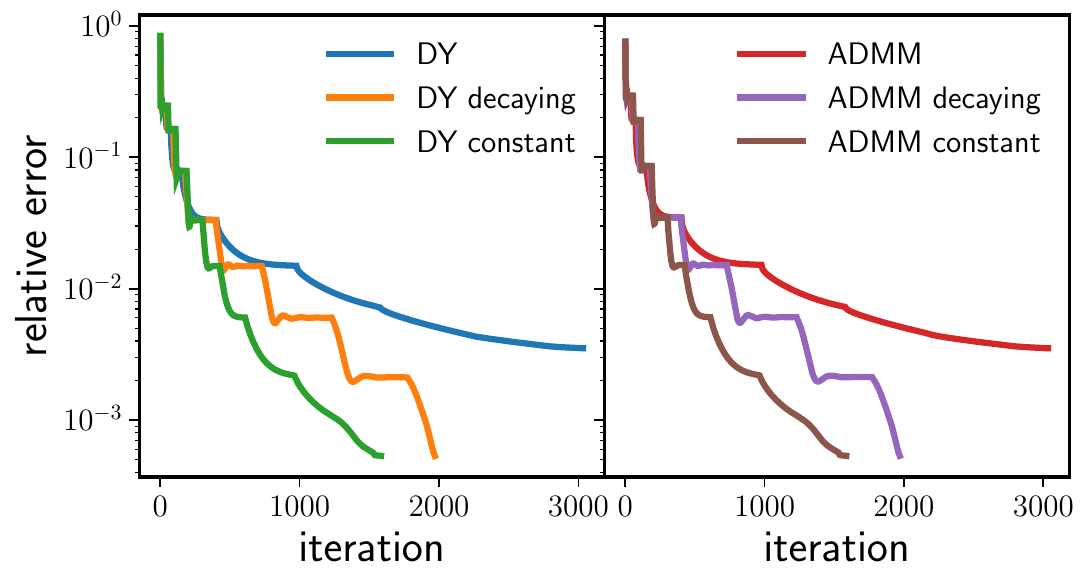}
\put(-480,5){\emph{(a)}}
\put(-240,5){\emph{(b)}}
\caption{\emph{(a)} Convergence rates of different methods when reconstructing
the observed image in Fig.~\ref{boats_fig}b by
solving \eqref{matcomp}. An example
of the recovered image is in Fig.~\ref{boats_fig}c.
\emph{(b)} Convergence rates for the same problem
but with annealing
on $\alpha$. An example of the recovered image with the accelerated
methods is shown in Fig.~\ref{boats_fig}d---the nonaccelerated
methods were unable to achieve such a small error.
\label{boat_convergence2}
\label{boat_convergence1}
}
\end{figure}

By  solving the nonnegative matrix completion
problem \eqref{matcomp}---with $a=0$ and $b=1$---our goal is to recover
$M$ from $M_{\textnormal{obs}}$.
As before, by running the previous algorithms a single time we
recover the image shown in Fig.~\ref{boats_fig}c. The convergence rates
of different methods are shown in Fig.~\ref{boat_convergence1}---in this
case we choose $\alpha=1$,
$\eta = 0.1$ for constant damping \eqref{hbdamp},
$r=3$ for decaying damping \eqref{nagdamp}, step size
$h=1$, and tolerance $\epsilon=10^{-6}$ in \eqref{tol}.
Note how the accelerated variants achieve faster convergence. %In  this  example,
All methods recovered the correct rank.

To obtain
more accurate results, we again consider annealing $\alpha$
according to \eqref{anneal_alpha}---we choose $\delta = 0.25$ and
$\bar{\alpha} = 10^{-4}$, which was the smallest feasible
value we found through several trials.
An example of the recovered image---with the accelerated variants since
the base methods were unable to achieve such a low error---is shown in
Fig.~\ref{boats_fig}d. By a closer inspection one can verify that
Fig.~\ref{boats_fig}d is actually better than Fig.~\ref{boats_fig}c, e.g.,
has a better resolution.
The convergence rates of each method under this setting  are shown in
Fig.~\ref{boat_convergence2}.
In this example, both accelerated variants---i.e., constant and decaying
damping---were able to significantly improve over the base methods.
Since this is a harder problem compared to the case of Fig.~\ref{mc2},
the objective function should have less curvature thus underdamping
the system yield better results.

\section{Discussion and outlook} \label{sec:conclusion}

The main goal of this paper was to provide a unified perspective
on standard optimization methods from a physics standpoint---these
results are summarized in Fig.~\ref{overview}.
We showed that there is indeed a unifying principle
behind several well-known, but also new, optimization methods.
The most general situation can be described by an underdamped Langevin
equation, which accounts for accelerated methods  in a
stochastic setting.
The randomness comes from a stochastic approximation
to gradients or proximal operators---%
this has similarities to a
Gaussian perturbation from a heat bath in a standard Langevin dynamics.
When the minibatch is sufficiently large, i.e., the ``temperature''
$T \sim S^{-1} \to 0$,
we recover the deterministic setting that is described by a dissipative
form of Newton's equation.
In addition,  in the high friction limit  $\eta \to \infty$,
the dynamics becomes nonaccelerated and thus
described by an overdamped Langevin or
by  a first-order gradient flow, in the stochastic and deterministic regimes,
respectively.
All the optimization methods we discussed arise as discretizations
of such systems.
Thus, the fact that many seemingly unrelated optimization algorithms consist of
actual \emph{simulations} of the same dissipative physical system is quite
surprising given that most of these methods were originally
introduced from a completely independent approach, and without
any connection to physics.

Using ODE splitting ideas and implicit discretizations,
we introduced several accelerated and stochastic generalizations of the most
important proximal  algorithms in the literature.
Moreover, important gradient-based methods, such as
gradient descent,
heavy ball, and Nesterov,
also fit our framework; in these cases there is
no splitting of the objective  function
and they correspond to explicit discretizations
(see the Appendix).
Interestingly, the heavy ball method turns out to be
a structure-preserving (i.e., conformal  symplectic) integrator,
while Nesterov introduces a spurious dissipation; see Ref. \cite{Franca:2019}
for details.
Let us also point out that recently we generalized symplectic integrators
to general dissipative Hamiltonian systems \cite{Franca:2020}, which  may
offer a systematic approach to construct new gradient-based
optimization methods based on simulations of physical systems.

The results of this paper explain why some different optimization methods
may actually  behave quite similarly in practice;
the reason is because they are
numerical integrators to the same physical system
and have the same order of accuracy.
Therefore, at this level, it is not possible to distinguish between different
methods. However, the connections we established allow one to perform
\emph{backward error analysis}, which would result into modified ODEs
that capture how a particular discretization perturbs the original
system---we carried out such an analysis in the case of
heavy ball and Nesterov \cite{Franca:2019}.
It would thus be interesting to consider a similar analysis
to some of the proximal methods previously discussed. This
would bring refined insights and a better understanding of these methods.
Another interesting question concerns the choice of damping.
We considered constant \eqref{hbdamp} and
decaying \eqref{nagdamp} damping
since these cases are well-understood, e.g., we know some convergence rates
and exact solutions in the quadratic case, as
described in Sec.~\ref{sec:rates}.
However, other choices are possible and
establishing the ``optimal'' damping strategy for a given
problem class is an interesting and nontrivial problem.

It is worth mentioning that
the connections between stochastic optimization and Langevin dynamics
offers an interesting opportunity for bringing
techniques from perturbation theory and nonequilibrium statistical
mechanics into machine learning.
For instance, any Fokker-Planck equation
is equivalent to a path integral through the Martin-Siggia-Rose
formalism \cite{ZinnJustin}; this is
a good starting point to apply field theory techniques to such problems.
However, the optimization methods we considered may
be applied to numerical problems in physics.
%For instance, Langevin dynamics simulation
%can be used to estimate path integrals and
%partition functions, thus , it should be possible to adapt
%these proximal optimization methods---by adding extra Gaussian noise,
%annealing schedule on the temperature, parallel tempering, etc.---to
%approach such types of problems that often arise statistical field
%theory. Perhaps, this may offer a cheaper
%alternative  compared to traditional Monte Carlo methods.
For instance, there is a close connection
between statistical
mechanics and combinatorial optimization problems, such as finding
ground states of spin glasses \cite{Parisi,Nishimori}.
The proximal methods considered in this paper
can be adapted to tackle combinatorial
problems such as max-cut and
graph partitioning, which are somehow equivalent.
Notably, in the context of machine learning and signal processing,
proximal  methods have been applied to large-scale
problems with  great success  and proved to be quite efficient.
It seems interesting
to consider whether such an approach can lead to
scalable and cheaper methods for disordered systems compared
to standard Monte Carlo approaches.

Finally, it is thus clear that standard optimization
methods are nothing but simulations of classical dissipative systems.
A natural question concerns extending this perspective
to the quantum level.
Can we construct quantum optimization algorithms
by simulating dissipative quantum systems?
Open quantum systems is still a controversial topic, but they may
offer an alternative to quantum annealing, whose discretizations may
yield general quantum optimization algorithms.
Thus, in the same way that a classical computer is employed to
implement standard---i.e., classical---optimization algorithms,
or equivalently  to simulate
a classical  dissipative  systems,
perhaps a quantum computer may be necessary to
simulate a dissipative quantum  system \cite{Feynman:1982}.

\bigskip

\subsubsection*{Acknowledgments}
\vspace{-.5em}
GF thanks Michael I. Jordan for helpful discussions and
Patrick Johnstone for comments on an earlier draft of this paper.
We also thank the anonymous referees for insightful comments that helped
to improve the paper.
This work was supported by grants
ARO MURI W911NF-17-1-0304, NSF 2031985, and NSF 1934931.

\bigskip

\appendix

\section[Stochastic gradient descent, heavy ball, and Nesterov]{Stochastic gradient descent, heavy ball, \\ and Nesterov}
\label{sec:gradient_based}

Since stochastic versions of gradient descent, heavy ball, and Nesterov's
method are widely used in machine learning, here we show how these methods
arise from the previous physical systems.

Let us start with gradient descent.
An explicit Euler-Maruyama discretization of the Langevin equation
\eqref{overLang} yields
\begin{equation}
\begin{split}
x_{k+1} &= x_k - h \nabla \F(x_k) + \dfrac{h}{\sqrt{S}} C(x_k) \epsilon_k \\
&= x_k - h \widetilde{\nabla \F}(x_k),
\end{split}
\end{equation}
where in the last passage we used \eqref{disc_noise}
and \eqref{sgrad_noise}.
This is nothing but the well-known \emph{stochastic gradient descent} (SGD)
method, widely used in training neural networks.
When $S\to\infty$, i.e. the stochastic gradient becomes
the true gradient, $\widetilde{\nabla \F} \to \nabla \F$,
we recover the deterministic gradient descent.

Let us now consider Nesterov's method.  Writing the underdamped Langevin
\eqref{underLang} as
\begin{equation}
\ddot{x} = - \nabla \F(x) - \eta(t) \dot{x} + \sqrt{\dfrac{h}{S}} C(x) \dot{W},
\end{equation}
where this should be understood in the It\^ o sense, we
discretize it with the help of \eqref{eulersec}, \eqref{disc_noise} and
\eqref{sgrad_noise} to obtain
\begin{equation}
\begin{split}
\dfrac{x_{k+1} - \hat{x}_k}{h^2} &= -\nabla \F(\hat{x}_k) +
\dfrac{1}{\sqrt{S}} C(\hat{x}_k) \epsilon_k  \\
&= - \widetilde{\nabla \F}(\hat{x_k}).
\end{split}
\end{equation}
Recalling the definition \eqref{xhat}, and redefining the step size as
usual, $h^2 \mapsto h$, we obtain
\begin{subequations}
\begin{align}
x_{k+1} &= \hat{x}_k - h \widetilde{\nabla \F}(\hat{x}_k), \\
\hat{x}_{k+1} &= x_{k+1} + \gamma_{k+1}(x_{k+1} - x_k).
\end{align}
\end{subequations}
This is a stochastic version of Nesterov's method---the original method is
recovered in the deterministic case where $S \to \infty$.

Let us now consider Polyak's heavy ball method.
From a physics point of view this case is more interesting since
the deterministic method was recently shown by us \cite{Franca:2020} to
be a conformal symplectic integrator.
Consider the underdamped Langevin equation \eqref{underLang}.
We write the second equation as
\begin{equation} \label{hb1}
\dfrac{d}{dt}\left( e^{\theta(t)}  p \right) = - e^{\theta(t)} \nabla \F(x)
+ e^{\theta(t)} \sqrt{\dfrac{h}{S}} C(x) \dot{W},
\end{equation}
where $\dot{\theta}(t) \equiv \eta(t)$.
Integrating this from $t_k$ to $t_{k+1} = t_k + h$, and keeping terms
up to $\bigO(h^2)$, we get
\begin{equation} \label{hb.u1}
\begin{split}
p_{k+1} &= e^{-h  \eta_k} p_k - h \nabla \F(x_k) + \dfrac{h}{\sqrt{S}} C(x_k)
\epsilon_k \\ %+ \bigO(h^2) \\
&= e^{-h \eta_k } p_k - h  \widetilde{\nabla \F}(x_k). %+ \bigO(h^2).
\end{split}
\end{equation}
The  first equation in \eqref{underLang} gives
\begin{equation} \label{hb.u2}
x_{k+1} = x_{k} + h p_{k+1} + \bigO(h^2).
\end{equation}
Now we define the following variables:
\begin{equation}
\mu_k \equiv e^{- h \eta_k}, \qquad v_k \equiv h p_k, \qquad  h^2 \mapsto h.
\end{equation}
This allows us to write \eqref{hb.u1} and \eqref{hb.u2} as
\begin{subequations}
\begin{align}
v_{k+1} &= \mu_k v_k - h \widetilde{ \nabla \F}(x_k), \\
x_{k+1} &= x_k + v_{k+1} .
\end{align}
\end{subequations}
This is precisely a stochastic version of the heavy ball method, often
referred to as
momentum method, or SGD with momentum, in deep learning literature.
Usually this method
is used with a constant $\eta_k = \eta$, in which case
$\mu \equiv e^{- h \eta} \in (0,1]$ is known as the ``momentum factor.''
When $S \to \infty$, i.e. the stochastic gradient becomes the true gradient,
the above recovers the deterministic heavy ball method.
Therefore, once again, the scheme of Fig.~\ref{overview} also apply
to gradient-based methods. Compared to the proximal methods previously
discussed, the difference in these cases is that the gradient is always
computed explicitly as opposed to implicitly.

\bibliography{biblio.bib}

\end{document}